\documentclass[final,5p,times,twocolumn]{elsarticle}

\usepackage{amssymb}

\usepackage{amsmath} 
\usepackage{mathtools} 	 	% :=  \coloneqq
\usepackage{algorithm2e}	% Pseudo Code
\usepackage{pgfplots}		% Plots in Latex
\usepgfplotslibrary{patchplots}	% 3D Plot
\usepackage{graphicx}
\usepackage{psfrag}
\usepackage{import}

\usepackage{textcomp}
\usepackage{siunitx}

\usepackage[utf8]{inputenc}
\usepackage{tikz}
\usetikzlibrary{trees}
\usepackage{amsthm}

\DeclarePairedDelimiter\floor{\lfloor}{\rfloor}

\journal{Computational Mechanics}

\newlength\figureheight
\newlength\figurewidth
\newcommand{\argmin}{\operatornamewithlimits{argmin}}

\newcommand{\dis}{\displaystyle}
\newcommand{\be}{\begin{equation}}
\newcommand{\ee}{\end{equation}}
\newcommand{\bea}{\begin{eqnarray}}
\newcommand{\eea}{\end{eqnarray}}
\newcommand{\balign}{\begin{align}}
\newcommand{\ea}{\end{array}}

\def\endbox{
            \leavevmode\kern-.25em\raise-.15em
            \vbox{\hrule height 0.8truept
            \hbox{\vrule width  0.8truept
            \vbox{\kern 4truept
            \hbox{\kern 1.5truept
                {}\kern 2.5truept}
                  \kern 3.0truept}
            \vrule width  0.8truept}
            \hrule height 0.8truept}
            }
\def\pmb#1{\setbox0=\hbox{#1}%
\kern-.025em\copy0\kern-\wd0
\kern.05em\copy0\kern-\wd0
\kern-.025em\raise.0433em\box0 }

\newcommand{\bD}{\mbox{\boldmath{$D$}}}

\newcommand{\bG}{\mbox{\boldmath{$G$}}}

\newcommand{\bK}{\mbox{\boldmath{$K$}}}

\newcommand{\bM}{\mbox{\boldmath{$M$}}}

\newcommand{\bP}{\mbox{\boldmath{$P$}}}

\newcommand{\bR}{\mbox{\boldmath{$R$}}}

\newcommand{\bU}{\mbox{\boldmath{$U$}}}
\newcommand{\bV}{\mbox{\boldmath{$V$}}}
\newcommand{\bW}{\mbox{\boldmath{$W$}}}

\newcommand{\bzero}{\mbox{$\bf 0$}}

\newcommand{\blambda}{\mbox{\boldmath{$\lambda$}}}

\newcommand{\bSigma}{\mbox{\boldmath{$\Sigma$}}}

\newcommand{\bPhi}{\mbox{\boldmath{$\Phi$}}}

\begin{document}

\begin{frontmatter}



% \bibliographystyle{unsrtnat}
\title{Proper orthogonal decomposition (POD) combined with hierarchical tensor approximation (HTA) in the context of uncertain parameters}


\author{Steffen Kastian, Dieter Moser, Stefanie Reese, Lars Grasedyck}

\address{RWTH Aachen, Germany}

\begin{abstract}
The evaluation of robustness and reliability of realistic structures in the presence of uncertainty 
involves costly numerical simulations with a very high number of evaluations. 
This motivates model order reduction techniques like the proper orthogonal decomposition. 
When only a few quantities are of interest an approximative mapping 
from the high-dimensional parameter space onto each quantity of interest is sufficient. 
Appropriate methods for this task are for instance the polynomial chaos expansion or 
low-rank tensor approximations. 

In this work we focus on a non-linear neo-hookean deformation problem
with the maximal deformation as our quantity of interest. 
POD and adaptive POD models of this problem are constructed
and compared with respect to approximation quality and construction cost.
Additionally, the adapative proper orthogonal decomposition (APOD) 
is introduced and compared to the regular POD.
Building upon that, several hierarchical Tucker approximations (HTAs) are constructed from
the reduced and unreduced models. 
A simple Monte Carlo method in combination with HTA and (A)POD is used to estimate 
the mean and variance of our quantity of interest.
Furthermore, the HTA of the unreduced model is employed to find feasible snapshots for (A)POD.

\end{abstract}


\end{frontmatter}



\section{Introduction}

An engineering task like optimization of structures, materials or manufacturing processes 
as well as the evaluation of robustness or failure probability of a system frequently relies on processing vast amounts of data.
In many cases, these data is generated from simulations, e.g. FEM, since real experiments are quite expensive and time consuming.
In order to obtain accurate results a huge amount of degrees of freedom (dof) are needed, 
which leads to time consuming simulations.

In contrast to real experiments, there are many approaches to reduce the time needed for a simulation.
A classical approach is the parallelization of the simulation. 
The success of this approach strongly depends on the scalability of the problem.
If e.g a Monte Carlo (MC) simulation is performed, each sample may be computed in parallel.
Therefore, the scalability is one. 
% Meaning that, doubling the number of computing units halfs the computation time.
On the other side of the spectrum one finds e.g. time stepping methods, 
where each time step fully depends on the preceeding time step. 
Without special methods this computation is not scalable. 
% Doubling the number of computing units does not reduce the computation time.
Assuming a high scalability, the next barrier is the cost of operating millions of computing units.
With around $1.7 MW$ power usage the JUQUEEN supercomputer at the JSC provided $458.752$ cores.
Depending on the problem this is not sufficient. Assume a system which is influenced by $10$ parameters, 
for each parameter we have $10$ options and the evaluation of the system for one parameter configuration takes $10$ hours.
Then, computing all paramater configurations would take around $24.88$ years on the JUQUEEN, 
more than four times of the life span of the machine itself.
This example illustrates the ubiquitous ``curse of dimensionality'' \cite{bellman2013dynamic}. 
Thus eventually, using parallelization on supercomputers reaches its limits and other approaches have to be considered.
One approach to tackle this curse is to reduce the computation time for one evaluation. 

Reducing the computation time for one evaluation is achieved by model order reduction (MOR) techniques 
like the proper orthogonal decomposition~(POD)~\cite{sachs2010pod,kunisch2008proper,kastian2018adaptive}, the proper generalized decomposition~(PGD)~\cite{chinesta2011short}, reduced basis~(RB)~\cite{eftang2011hp} and the empirical interpolation method~(EIM)~\cite{radermacher2016pod}.
PGD introduces different parameters as extra-coordinates~\cite{chinesta2011short}. 
% Furthermore, it is neccessary to set up the PGD for each problem separately. 
% These are both drawbacks in comparison to POD.
More general are projection based MOR techniques like 
POD~\cite{berkooz1993proper}, load-dependent Ritz methods~\cite{Krysl2001} and modal basis reduction~\cite{nickell1976nonlinear}, 
which differ in the projection matrix used.
The empirical interpolation method is an extension of POD. 
EIM  splits the operators in a linear and non-linear part. 
The linear part is treated by POD, the non-linear part by an additional subspace projection~\cite{barrault2004empirical}.
%Modifying the POD in a sense of reducing the complexity of evaluating the nonlinear term of the reduced model leads to the empirical interpolation method \cite{barrault2004empirical}. 

%\textcolor{red}{Beginn Ideen:}
%K	POD - Projection based
%K	PGD - 
%reduced basis
%K	empirical interpolation
%modal reduction
%balance truncation
%Krylov methods
%moment matching
%\textcolor{red}{Ende Ideen.}
%\textcolor{red}{An Steffen, versuch mal in 4-5 Sätzen einen überblick über diese oder noch mehr Methoden dieser Art zu geben! Mit Referenzen.}

Even if one finds a sufficiently accurate reduced model and thereby reduces the computational time substantially,
the curse of dimensionality remains. 
For $10$ parameters with $20$ options and
an evaluation time of only $1$ millisecond, the parallelized evaluation of all parameters would still take $6912$ years.
The only way to reduce this time is to reduce the number of evaluations. 
How this can be achieved, depends on the question posed.

In many situations, a simulation can be reduced to one or several characteristic values. 
In a production chain this could be the number of produced units per day.
For the design of a CPU this could be the power usage per operation, and for
the construction of a bridge this could be the maximal stress over the day.
Often, this characteristic value is denoted as quantity of interest (QoI).
Then, a task in practical design might be to find a parameter configuration which results into maximization or minimization of this QoI.
Uncertainty quantification includes the computation of statistical moments or level sets of the QoI.

When a mapping from paramater space onto the QoI is given in a proper form, 
it becomes feasible to perform these tasks without the evaluation of every paramater combination.
Let $p_1,\ldots p_d$ be the $d$ parameters and imagine that the QoI $\varphi$ is well approximated by
\begin{align*}
	\tilde{\varphi}\left( p_1,\ldots,p_d \right) = \prod_{k=1}^{d}f_k(p_k).
\end{align*}
In this fully separated form, the computation of the absolute maximum is reduced to
finding the absolute maximum of the $d$ functions $\left\{ f_k \right\}_{k\in \left\{ 1,\ldots,d \right\}}$.
Similarly, under the assumption that the parameters are statistically independently distributed
the computation of the mean reduces to the multiplication of the $d$ means of $\left\{ f_k \right\}_{k\in \left\{ 1,\ldots,d \right\}}$.
Obviously, in practice it is not always possible to find this particular form of separability.
Nonetheless, many other forms of separability are found in various tensor formats, e.g., 
the Canonical Polyadic format~\cite{hitchcock1927expression}, 
the Tucker format~\cite{de2000multilinear}, 
the Tensor Train format~\cite{oseledets2011tensor} or the Hierarchical Tucker format~\cite{hackbusch2009new,grasedyck2011introduction,grasedyck2010hierarchical}. 
This class is denoted as low-rank approximations.
In \cite{kolda2009tensor,hackbusch2012tensor} very good overviews are found.
These formats are actively researched and were already applied successfully to 
mitigate the curse of dimensionality in many situations, examples are~\cite{ballani2015hierarchical,litvinenko2012uncertainty,espig2014efficient,matthies2016parameter,dolgov2017hybrid,corveleyn2017computation}.

In this work methods from model order reduction and low rank approximations are investigated for a Neo Hookean model problem, which is presented in Sec.~\ref{sec:cube_under_compression}. 
From the class of MOR techniques, the novel adaptive proper orthogonal decomposition APOD is introduced in Sec.~\ref{sec:APOD}. 
The advantages of APOD over POD are worked out in Sec.~\ref{sec:APOD_result}. 
Especially in the case of additional geometric parameters the advantage is clearly shown in the numerical results.
These numerical results also motivate the usage of low-rank approximations.
From the latter class, the Hierachical Tucker approximation (HTA) is briefly introduced in Sec.~\ref{sec:ht_theory}.
The HTA is applied to the model problem and the approximation quality is investigated in Sec.~\ref{sec:HTA_cube}. 
A first idea to find synergies between HTA and APOD is established in Sec.~\ref{sec:snapshot_search}.
In that section a low approximation of the POD error is constructed, 
in order to find paramater configurations and thereby snapshots, which improve the POD model. 
The uncertainty quantification tasks performed in this work are motivated in Sec.~\ref{sec:UQ_for_cube}.
%In Sec.\ref{sec:UQ_with_HTA} we describe how to use the HTA for these tasks.
Finally, in Sec.~\ref{sec:comparison_HTA_APOD} these tasks are performed with POD, APOD and HTA and thereby compare our methods.
The paper is concluded with an outlook in Sec.~\ref{sec:outlook}.

\section{Model problem}
\label{sec:cube_under_compression}

\begin{figure}
%      \psfrag{x}{$x$}
%      \psfrag{y}{$y$}
%      \psfrag{z}{$z$}
%      \psfrag{P}{$P$}
%      \psfrag{2l1}[c]{$2l_1$}
%      \psfrag{2l2}{$2l_2$}
%      \psfrag{l3}{$l_3$}
     \includegraphics[width=7cm]{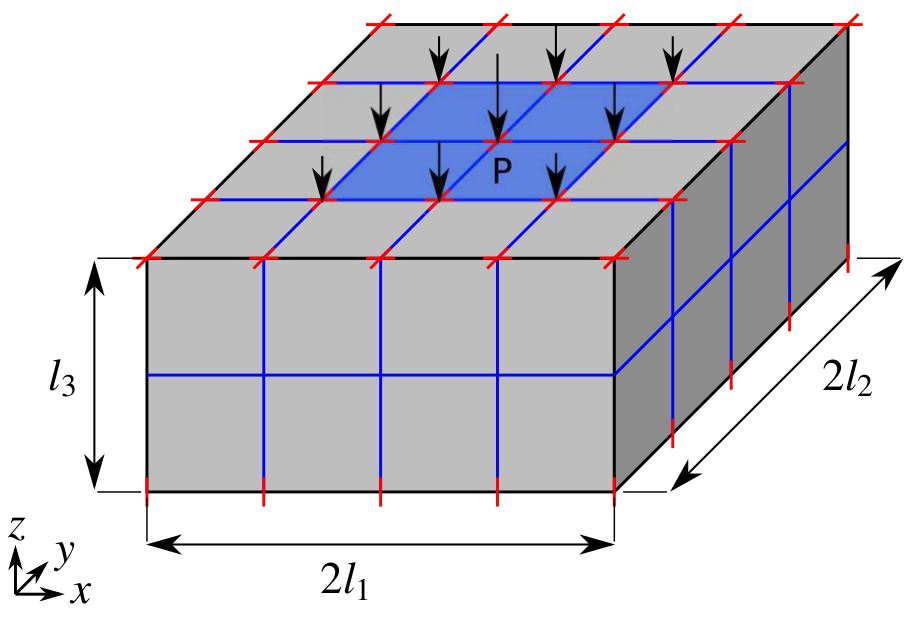}   
   \caption{Cube under compression}
   \label{fig:cube_under_compression}
\end{figure}

The model problem investigated here is a cube made of Neo-Hookean material undergoing large deformation. 
Geometry and boundary conditions are shown in Fig.~\ref{fig:cube_under_compression} (c.f. \cite{reese2000}).
% The model problem we present in this paper is the cube under compression with a Neo-Hookean material law and large deformation 
% with the geometry and loading conditions as depicted in Fig.~\ref{fig:cube_under_compression}. 
Using the symmetry of the problem only a quarter of the cube is considered.
In Tab.~\ref{tab:cube_under_compression} the parameters of our main setup are found.

% \begin{figure}[ht]
%  \psfrag{sym}[cc][cc][1]{sym.}
%  \psfrag{sym2}[cc][cc][1]{ sym.}
%  \psfrag{x1}{$x_1$}
%  \psfrag{x2}{$x_2$}
%  \psfrag{x3}{$x_3$}
%  \psfrag{a}[tc][cc][1][-10]{$l_1$}
%  \psfrag{b}[tc][cc][1][15]{$l_2$}
%  \psfrag{c}[bc][cc]{$l_3$}
%  \includegraphics[width=7cm]{pictures/block_under_compr/block_under_compr.eps} %block_under_compr
%  \caption{Cube under compression}
%  \label{fig:cube_under_compression}
% \end{figure}

% In a first step a change of the Young's modulus is considered for different number of modes. 
% Therefore, all other parameters are the same as the precalculation found in Tab.~\ref{tab:cube_under_compression}.

%\begin{table}[ht]
% \caption{Material parameters main precalculation}
% \label{tab:cube_under_compression}
% \begin{tabular}{lll}
%  Material 	& Parameters 			& Pressure		\\
%  Neo-Hookean	& $E=90 N /mm^2$		& $-0,320N/mm^2$ 	\\
%		& $\nu=0.4999$			&			\\
%		& $l_1=l_2=l_3=1.0mm$		&			\\
% \end{tabular}
%\end{table}

\begin{table}[ht]
 \begin{tabular}{llll}
	 Young's module & $E$ &$=90 $ & $ \left[ \SI{}{\N / \mm\squared} \right] $ \\[1ex]
	 Poisson's ratio& $\nu$ &$=0.4999 $ & \\[1ex]
	 geometric dimensions& $l_1,l_2,l_3$ &$= 1 $ & $ \left[ \SI{}{\mm} \right] $ \\[1ex]
	 Pressure & $F$ &$= -320 $ & $ \left[ \SI{}{\N / \mm\squared} \right] $ \\
 \end{tabular}
 \caption{Parameter configuration for the main precalculation.}
 \label{tab:cube_under_compression}
\end{table}

To ensure a stable computation, the pressure is applied incrementally in time. 
100 time steps are used until full loading is achieved
%\textcolor{red}{Wieviele Zeitschritte und wie lang sind diese? / Nicht konstant für das Paper. Meist 100\% loading => 100 time steps}
%\textcolor{red}{
For the computations itself the finite element analyse program \textsc{FEAP}\cite{taylor2000feap} is used.
%It is therefore possible to add all kinds of extensions, like MOR and element subroutines like damage modeling.$(*1)$ }
%\textcolor{red}{ very good program with very wise parameters and nice elements on very fast computers. Okay like this?}
%A accurate computation with a fine discretized mesh with FE takes a long time. 
%When thousands of FEAP runs are needed 
%for optimization algorithms or Monte Carlo simulations, this becomes infeasible.
%The problem aggravates, when the geometry becomes more challenging
%and therefore the number of elements increases. 
%In both cases we are strongly interested in reducing the complexity of the problem.

% Vllt.: We use the output of these FE computations as our unreduced model.

\subsection{UQ for the Cube}
\label{sec:UQ_for_cube}

Since the paper focuses on surrogate models for uncertainty quantification, cf.~\cite{smith2013uncertainty}, as main application,
our problem is extended, so that the parameters of our model problem lie in a range
instead of having an assigned fixed value. 
Additionally, the units are omitted and 
the parameter ranges are set to
\begin{align}
	E \in \left[ 150, 2000 \right] \mbox{ and }
	l_k \in \left[ 0.95, 1.05 \right],\, k\in\left\{ 1,2,3 \right\}.
	\label{eq:parameter_ranges}
\end{align}
%\textcolor{red}{Teilweise E-Modul von 90 bis 5000? Kürzen auf bis 2000?}
For the remainder of the paper we will keep the boundaries for the geometric parameters constant,
the Young's modulus and the loading will change to save some computation time.
To ease up notation in the following chapters, we also introduce the parameter space
\begin{align}
	\mathcal{S} = \left[ 0.95,1.05 \right]^{3} \times \left[ 150,1000 \right]
	\label{eq:parameter_space}
\end{align}
and take the liberty to enumerate our parameters
\begin{align}
	p_k = l_k, \mbox{ for } k\in\left\{ 1,2,3 \right\} \mbox{ and } p_4 = E.
	\label{eq:parameter_names}
\end{align}

%\textcolor{red}{Nicht sicher ob wir das so machen}.

The problem of uncertainty quantification eases up, 
when only one scalar value is of interest. 
For our model problem the maximal deformation is chosen, which is 
measured in percentage of the initial height of the cube.
This quantity of interest is denoted as $\varphi$.
This leads us to the main task of this paper, 
which is to provide a mapping 
\begin{align}
	\varphi : \mathcal{S} \to \mathbb{R}, p \mapsto \varphi(p) .
	\label{eq:main_mapping}
\end{align}%
Once a mapping is found, it becomes feasible to estimate 
the first statistical moments of the quantity of interest,
namely the mean and the variance, by means of a Monte Carlo method. 
Each parameter has a probability distribution, e.g. a uniform distribution over the range or
a truncated normal distribution. 
We generate a finite discrete set $\mathcal{X}$ of random parameter configurations
according to the respective distributions. 
Then, the sums
\begin{align}
	\mu = \frac{1}{| \mathcal{X} |}\sum_{p \in \mathcal{X}} \varphi(p) \mbox{ and } \sigma^{2} = \frac{1}{| \mathcal{X} |-1}\sum_{p \in \mathcal{X}} \left( \varphi(p) - \mu \right)^{2}
	\label{eq:mc_mean_var}
\end{align}%
estimate the mean and the variance.
Non-intrusive methods like multilevel Monte Carlo~\cite{giles2013multilevel}
or stochastic collocation~\cite{babuvska2007stochastic} 
methods are superior to this method 
regarding the number of samples needed for an accurate estimation.
But since the focus of this paper lies on the construction and investigation of reduction methods for the mapping \eqref{eq:main_mapping} 
the standard Monte Carlo method is used throughout this paper. 
In the following chapters, the reduced models which provide this mapping are studied in detail. 
Then, then the resulting methods are employed for uncertainty quantification in Section \ref{sec:comparison_HTA_APOD}.
%\textcolor{red}{If the reduced models are suited for the Monte Carlo method, then it is to assume that this yields improvements when multilevel Monte Carlo or stochastic quadrature is used. $(*2)$ }
%\textcolor{red}{Zitate Zitate Zitate}

\section{Model reduction for the Neo Hooke equation}
\label{sec:MOR_for_NHE}
% \textcolor{red}{1) MOR: (Oder eher Darstellung des eigentlichen Problems)}
% Wie macht man das grob
Model order reduction (MOR) in the context of FE means that the order of the model is reduced in order to save CPU time and memory \cite{radermacher2014model}.
For many applications this paves the way for indepth investigations of the problem at hand, 
because irrelevant parts of the model can be significantly reduced.

\subsection{The problem}
% Projection based model order reduction (MOR) in the context of FE means that the high dimensional equation system is projected to a equation system of smaller dimension \cite{krysl2001dimensional}. This smaller equation system is solved and the solution can be projected back to the high dimensional space.

% General equation System 
The discrete nonlinear equation system is obtained from a displacement-based finite element formulation  %We begin with the equation system 
\be \label{eq:generalnonlinear}
 \bM \ddot{\bU}(t) +\bR(\bU(t))=\bP (t),
\ee
where for static examples the $n \times 1$ nodal acceleration vector $\ddot{\bU}(t)$ can be neglected. 
The $n \times n$ matrix $\bM$ represents the mass matrix and the $n \times 1$ vector $\bU (t)$ the vector of nodal displacements. 
The external load is represented by the time-dependent $n \times 1$ vector $\bP (t)$.
With this assumption and Eq. \eqref{eq:generalnonlinear} the residual vector $\bG(\bU(t))$ is defined as 
\be \label{eq:definitionG}
  \bG (\bU(t)) \coloneqq\bR(\bU(t)) - \bP(t)= \bzero.
\ee
% General Newton-Raphson
 
%  It has been shown in previous researches by \citet{radermacher2013comparison} that the POD method is an appropriate possibility to get cheaper simulations with a very good accuracy. 
The dependence on the time $t$ is from now on omitted in the notation.
Commonly, the Newton-Raphson method is used to solve the nonlinear Equation \eqref{eq:definitionG}.
% Using the Taylor expansion of $\bG (\bU_{i+1}^j)$ yields the well known iterative solution algorithm
% \begin{align}\label{eq:iterative_solution_alg}
%   \bG (\bU_{i+1}^j) \approx \bG(\bU_i^j)+\bK_T(\bU_i^j) \Delta \bU_{i+1}^j = \bzero \nonumber \\
%    \bU_{i+1}^j = \bU_i^j+\Delta \bU_{i+1}^j \nonumber \\
%    ||\bG(\bU_{i+1}^j)|| \le  \tol \nonumber \\
%    i \leftarrow i+1 ,
% \end{align}
% where $$\bK_T(\bU_i^j) = \left. \frac{\partial \bG\left (\bU\right )}{\partial \bU}\right|_{\bU_i^j}$$ is the tangential stiffness matrix with the dimension $n \times n$. 
% In this notation $j$ refers to the $j$th time step and $i$ to the $i$th Newton iteration step. 
% Projection matrix - Model order reduction
\subsection{Projection based MOR}
% The number of unknowns $n$ which is equivalent to the number of degrees of freedom can be very high which leads to long computational times. To reduce this problem one can use model order reduction method (MOR) techniques. 
One special kind of MOR is the projection based MOR, where the high dimensional equation system is projected to an equation system of smaller dimension $m$ \cite{krysl2001dimensional}. 
This smaller equation system is solved and the solution is then projected back to the high dimensional space. 
%Which means that the high dimensional equation system is represented by a linear combination of the vectors of the projection matrix also known as basis function or modes. 
%This works good for problems where most dofs stand in relation to each other during the different timesteps of a simulation. 
Especially suitable is this procedure for problems, where the difference in displacement between time steps is well approximated by a linear combination of the basis $\bW_1, \ldots ,\bW_m$. 

The dimension $m$ determines the accuracy and computational effort of the POD. 
With smaller $m$ the computational effort is reduced in exchange for accuracy.
Usually, $m$ can be chosen to be much smaller than $n$, without substantial loss of accuracy. 

Let $\bPhi$ be the $n \times m$ projection matrix which serves for the projection from the $n$ dimensional space onto the $m$ dimensional space.
Furthermore, let $\bU_{\mathrm{red}}^j$ be the reduced $m \times 1$ displacement vector and 
\be \label{eq:Ured}
 \bU^j \approx \bar{\bU}^j = \bPhi^j \; \bU_{red}^j
\ee
be the unreduced $n\times 1$ displacement vector 
of the $j$th time step. 
Additionally, the tangential stiffness matrix $\bK$ can be projected as
\be \label{eq:Kred}
\bK_T(\bU^j) = \bPhi^j \bK_T(\bU_{red}^j) \bPhi^{j~T}
\ee
where $\bK_T(\bU_{red}^j)$ is the reduced tangential stiffness matrix. 
Applying Eq. \eqref{eq:Ured} and Eq. \eqref{eq:Kred} to the Taylor-Expansion Eq. \eqref{eq:iterative_solution_alg} 
the dimension is reduced from $n$ to $m$:
\be
 \underbrace{\bPhi^{j~T} \bG( \bPhi^j \bU_{red,i}^j)}_{\dis{\bG_{red}}} + \underbrace{\bPhi^{j~T}\bK_T(\bPhi^j\bU_{red,i}^j)\bPhi^j}_{\dis{\bK_{t,red}}} \Delta \bU_{red,i+1}^j = \bzero.
\ee
% \textcolor{red}{(Nochmal Newton Raphson reduced?)}

% \textcolor{red}{(Referenzen?)}
The success of this procedure relies on the right choice of the $m$ dimensional space and thus
respectively on the proper construction of the projection matrix $\bPhi$.
%To project to a reduced equation system it is neccassary to construct a proper projection matrix $\bPhi$. 
Alternative methods to construct the projection matrix are the modal basis reduction method~\cite{nickell1976nonlinear} and the load-dependent Rith method~\cite{mollemans2005very,krysl2001dimensional}.
For the type of problems which are her of interest, the POD has been shown to be the most efficient~\cite{radermacher2013comparison}.

\subsection{Proper orthogonal decomposition (POD)}\label{sec:POD}

The proper orthogonal decomposition method (POD) is based on the singular value decomposition (SVD). % like the Karhunen-Loève expansion (KLE) or prinicipal component analysis (PCA). 
It is mainly used in fields like turbulent flow, stochastics, image or signal analysis~\cite{chatterjee2000introduction,holmes2012turbulence}. 
In the field of mechanics \cite{berkley2004real, kerschen2005method, lenaerts2001proper} different approaches to employ POD were presented. 

The POD is built upon ``snapshots'' which are the result of unreduced calculations.
Thus, it is important to perform the unreduced calculations with parameter configurations
which in the end represent the full behaviour of the system.
%In our case the displacements of different time steps of precalculations are saved as snapshots. 
%It is crucial to select proper parameter settings as precalculations. 
Each time step of our precalculation results into a displacement vector $\bU^j$ with the dimension $n \times 1$.
These displacement vectors are gathered in the $n \times l$ snapshot matrix 
$ \bD = \left[ \bU^1 ~ \bU^2 ~  ... ~ \bU^{l}\right] ~\in \mathbb{R}^{n \times l}.$
The SVD 
% Having these information a suitable orthonormal basis $\bPhi$ needs to be found. Computing a singular value decomposition of $\bD$ can solve this problem
$  \bD = \bW \bSigma \bV^T = \sum_{k=1}^{l} \sigma^k \bW^k \bV^{k^T}$
is used to generate a suitable orthonormal basis where $\bW = (\bW^1~\bW^2 ~...~ \bW^l) ~ \in \mathbb{R}^{n \times l}$ 
and $\bV = (\bV^1~\bV^2 ~...~\bV^n) ~ \in \mathbb{R}^{l \times n}$ are orthonormal matrices. 
$\bSigma$ is a diagonal matrix which contains the singular values ($\sigma^k$) which are always positive, ordered decreasingly. The projection matrix $\bPhi$ can be obtained as
\be
 \bPhi = \left( \bW^1 ~ \bW^2 ~ ... ~ \bW^m \right) ~~\in  \mathbb{R}^{n \times m}.
\ee
The projection matrix $\bPhi$ is used to solve the nonlinear Eq. \eqref{eq:definitionG} on a lower dimensional space.

\section{Adaptive proper orthogonal decomposition (APOD)}\label{sec:APOD}

% \textcolor{red}{Warum schlecht für Nichtlineare Probleme? Motivation?}
% Previous works \cite{kastian2018adaptive} shown that POD is not always suitable for non-linear problems.
When using POD, it is assumed that the solution lies in a linear space, spanned by the basis functions $\bW^1 \ldots \bW^m$. 
Therefore, each linear combination of these basis functions is a possible solution. 
In general, all snapshots are used to construct the projection matrix $\bPhi$.
However, in particular in the context of non-linear behaviour, not all snapshots might be representative of the current deformation behaviour.
This makes the POD inefficient, since more basis vectors than necessary are used~\cite{haasdonk2011training}.

It is imaginable that in each time step the respective solution would be matched better by a different but smaller basis. 
Thus, a simple principle to choose a proper subset of snapshots is introduced
which are then employed in the construction of such a matched basis. 
We call this method \textit{adaptive proper orthogonal decomposition}.

\subsection{Principle}
For geometrically nonlinear problems we suggest to introduce a projection matrix $\bar{\bPhi}$ which depends on the maximal current displacement $\bar{u}_{max}$.
Investigating geometric nonlinear problems we propose to introduce a projection matrix $\bar{\bPhi}$ depending on the maximal current deformation $ \bar{u}_{max} $, 
%For each timestep one get a new projection matrix which fits better to the current state. 
the bar represents the current time step and the maximum displacement is defined as
\be
 \bar{u}_{max} = \| \bar{\bU} \|_{\infty}.
\ee
%In the following we derive a projection matrix depending on $\bar{\bU}$ to make sure that the projection suits the current state. 
%Other possibilities like taking the absolute value of the maximum displacement of a node or to compare the angles between the current state and the snapshots could work too.
The snapshot vectors are sorted with respect to the maximum displacement
\be
 ||\bU^1||_{\infty} \le ||\bU^2||_{\infty} \le ... \le ||\bU^l||_{\infty}.
\ee
The current maximal displacement $u_{max}$ is compared with the maximal displacement 
of each of the snapshots $||\bU^1||_{\infty}~...~||\bU^l||_{\infty}$. 
The index of the snapshot, which fits the maximum displacement best is denoted as
\be
 b (\bar{u}_{max})=\underset{i\epsilon\{ 1,...,l\}}{\argmin}~\left(|(||\bU||_{\infty}-\bar{u}_{max})|\right).
\ee
The snapshots close to the current maximal displacement are chosen
and it is expected that these snapshots represent the behaviour of the current state well. 
To ensure that the snapshots around the index $b$ are still in the range of available snapshots the index of the first taken snapshot is defined as
\be
 a = \max{\left[ 0 ; \min \left(l-o;b- \floor*{\frac{o}{2}}\right)\right]}
\ee
which ensures that the snapshots exist. 
The number of considered snapshots is called  $o$. 
In \cite{kastian2018adaptive} it is shown that
\be
 o = ms
\ee
is a good choice for $o$, 
where $s$ is the number of precalculations and $m$ the number of modes. 
In the case that $$l-o<b-\floor*{\frac{o}{2}}$$ available snapshots 
do not lead to good results, we recommend to generate snapshots with larger deformations. 
After determining $a$ one can construct a suitable snapshot matrix  
\be
 \bar{\bD} = \left( \bU^{a+1} ~ \bU^{a+2} ~ ...~\bU^{a+o} \right) \in \mathbb{R}^{n \times m}
\ee
for the current timestep.
With a singular value decomposition of the current snapshot matrix 
\be
  \bar{\bD} = \bar{\blambda} \bar{\bSigma} \bar{\bV}^T = \sum_{k=1}^{l} \bar{\sigma}^k \bar{\blambda}^k \bar{\bV}^{k^T} 
\ee
the projection matrix $\bar{\bPhi}$ for the current timestep is created as
\be
 \bar{\bPhi} = \left( \bar{\blambda}^1 ~ \bar{\blambda}^2 ~ ... ~ \bar{\blambda}^m \right),
\ee
where $\bar{\bSigma}$ contains the singular values in decreasing order.
% \subsubsection{Adaptive number of modes}
% The number of modes affects both the accuracy and the calculation time. It can be chosen a priori or chosen with a criterion based on the fact that the accuracy depends on the decay of the singular values\cite{ammar2009coupling}. 
% \be
%  \frac{\bar{\sigma}_m}{\sum \bar{\sigma}} \le tol
% \ee
% Regarding the previous equations this has to be an iterative criterion. Means one starts with $m=3$ do all neccassary steps and check the criterion. If this criterion is fullfilled the projection matrix $\bar{\bPhi}$ yields accurate results. If the criterion is not fullfilled more modes are required which means
% \be
%   m \leftarrow m+1.
% \ee
\subsection{Numerical results}\label{sec:APOD_result}
% In section \ref{sec:APOD_result} the adaptive proper orthogonal decomposition (section \ref{sec:APOD}) is compared with the classical proper orthogonal decomposition (section \ref{sec:POD}) regarding accuracy and computational time. First we discuss the cube under compression for a variation of the Young's modulus $E$ with a fixed snapshot matrix $\bD$ and a fixed size of the cube. In a next step the possibility of uncertain parameters regarding the size of the cube are investigated. 
% In this paper we present the investigated cube under compression with a Neo-Hooke material law and large deformation. To study the accuracy and time efficiency of APOD we studied a geometry with boundary and loading conditions as depicted in Fig.~\ref{fig:cube_under_compression}. 
The parameters for the precalculation can be found in Tab.~\ref{tab:cube_under_compression}.
In a first step a change of the Young's modulus is considered for different number of modes. Therefore, all other parameters are the same as the precalculation found in Tab.~\ref{tab:cube_under_compression}.
In the following different models will be compared. % with respect to each other. 
To do so a relative error between two models X and Y is introduced as
\be \label{eq:Error}
	\epsilon_{ \mathrm{X,Y} } = \frac{ \varphi_{ \mathrm{X} }  -  \varphi_{ \mathrm{Y} } }{ \varphi_{ \mathrm{X} }  }.
\ee
The result of the FEAP simulation is denoted as $\varphi_{\mathrm{full}}$ where phi is defined in Section \ref{sec:cube_under_compression} 
and serves here as a placeholder for our quantities of interest.
The results of APOD and POD are denoted as $\varphi_{\mathrm{POD}}$ and $\varphi_{\mathrm{APOD}}$, respectively.
The abbreviation (A)POD is used, if both APOD and POD are referred to.
The error $\epsilon_{\mathrm{full,(A)POD}}$ is used to assess the accuracy of APOD and POD.

\begin{figure}[ht]
%  \begin{tikzpicture} 
%    \begin{semilogyaxis}[line width=0.2mm, xmin=100, xmax=2000, ymin=0.0000000001,extra x ticks={100}, legend style={at={(0.5,-0.22)}, anchor=north,legend columns=4}, xlabel={Young's modulus E  [$\SI{}{\N / \mm\squared}$] }, ylabel={$\epsilon_{\mathrm{full,(A)POD}}$} ] 		    %xmin = 90, xmax=5000
%      \addplot[color=blue, mark=x] table {Esnap_90-APOD3}; 
% 	\addlegendentry{APOD3} 
%      \addplot[color=blue, mark=*] table {Esnap_90-APOD5}; 
% 	\addlegendentry{APOD5} 
%      \addplot[color=blue, mark=square] table {Esnap_90-APOD7}; 
% 	\addlegendentry{APOD7} 
%      \addplot[color=blue, mark=triangle] table {Esnap_90-APOD10}; 
% 	\addlegendentry{APOD10} 
%      \addplot[color=blue, mark=o] table {Esnap_90-APOD20}; 
% 	\addlegendentry{APOD20} 
%      \addplot[color=red, mark=x] table {Esnap_90-POD3}; 
% 	\addlegendentry{POD3} 
%      \addplot[color=red, mark=*] table {Esnap_90-POD5}; 
% 	\addlegendentry{POD5} 
%      \addplot[color=red, mark=square] table {Esnap_90-POD7}; 
% 	\addlegendentry{POD7} 
%      \addplot[color=red, mark=triangle] table {Esnap_90-POD10}; 
% 	\addlegendentry{POD10} 
%      \addplot[color=red, mark=o] table {Esnap_90-POD20}; 
% 	\addlegendentry{POD20} 
% %       \legend{$APOD3$,$APOD5$,$APOD7$,$APOD10$,$APOD20$} 
%     \end{semilogyaxis} 
%    \end{tikzpicture}
     \includegraphics[width=7cm]{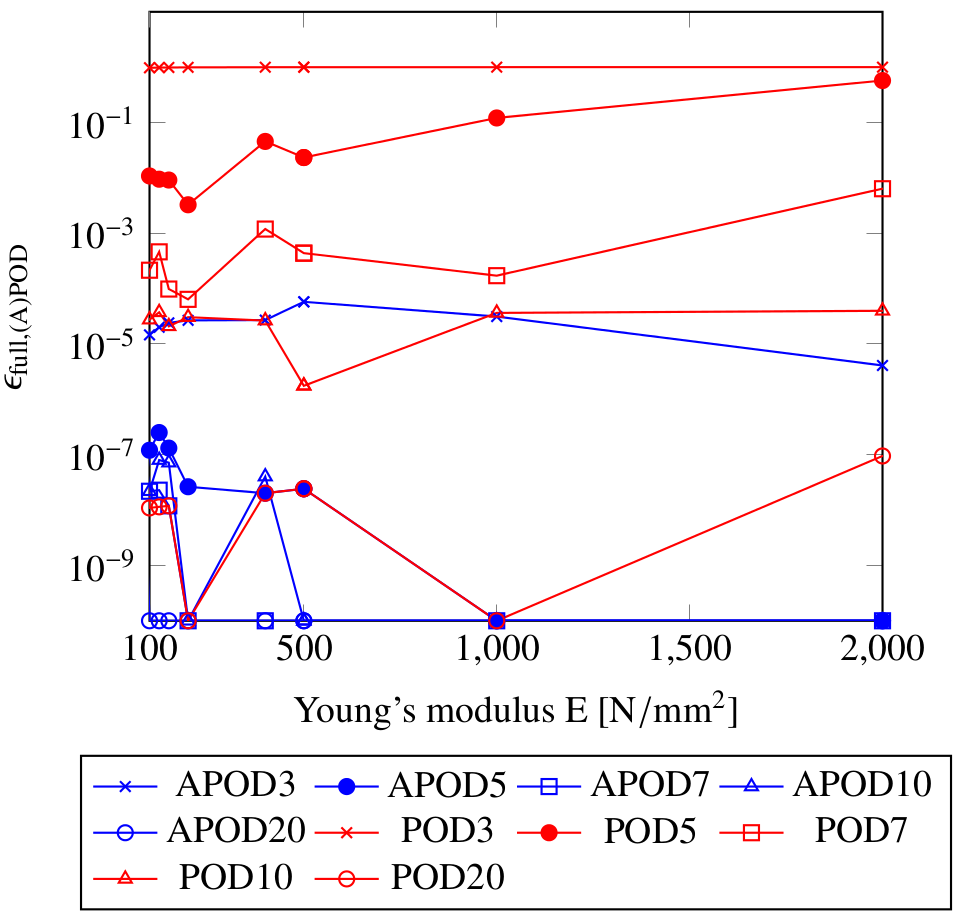}   
 \caption{Comparison of the accuracy of POD and APOD for a variation of the Young's modulus $E$.
 	  The number appending ``POD'' or ``APOD'' represents the number of used modes.}
 \label{fig:cube-E90-2000}
\end{figure}

In Fig.~\ref{fig:cube-E90-2000} POD and APOD are compared. 
APOD with $m=3$ modes is denoted by ``APOD3'', POD with $m=5$ modes is denoted by ``POD5''. 
The other configurations are denoted accordingly. 
In this example it can be seen that the new APOD method with the same number of modes is for all cases significantly more accurate than the POD method, since the snapshots are carefully selected to fit the current situation in every time step. 
Figure~\ref{fig:cube-E90-2000} also shows that the APOD method with $m=3$ modes yields approximately the same accuracy as POD with $m=10$ modes. 
APOD is more expensive due to the fact that for each time step a singular value decomposition is computed. 
On the other hand, each singular value decomposition is much cheaper since $\dim{(\bar{\bD})}\le\dim{(\bD)}$. 
The comparison of the computational time $t_{\mathrm{POD}}$ for POD and the computational time $t_{\mathrm{APOD}}$ for APOD shows that the cost are almost equal for the same number of modes, with  $t_{APOD}/t_{POD} \in \left[ 1.03 , 1.05 \right]$.
Comparing the computation times which yield approximately the same accuracy shows that with $t_{POD10}/t_{APOD3} \approx 1.63$ the POD requires 63\% more computation time.

%\textcolor{red}{Hier noch was zum Bild \ref{fig:l2_vs_qoi}. Legend gestrichelt L2}
% Plot L2 vs qoi
% \begin{figure}
% 	\centering
% 	\setlength\figureheight{0.45\textwidth*5/7}
% 	\setlength\figurewidth{0.45\textwidth}
% 	\InputIfFileExists{../HTA_choose_snapshots/uncertainty_1_precomp/1_precom.tikz}{}{\textbf{!! Missing graphics !!}}
% 	\caption{ Comparison of the $L2$ and QoI error }
% 	\label{fig:l2_vs_qoi}
% \end{figure}
% \begin{figure}
% 	\centering
% 	\setlength\figureheight{0.45\textwidth*5/7}
% 	\setlength\figurewidth{0.45\textwidth}
% 	\InputIfFileExists{../HTA_choose_snapshots/uncertainty_1_precomp/1_precom_qoi_rel.tikz}{}{\textbf{!! Missing graphics !!}}
% 	\caption{ Comparison of the relative QoI error }
% 	\label{fig:qoi_rel}
% \end{figure}
\begin{figure}
%  \begin{minipage}{0.2\textwidth}
	\centering
% 	\setlength\figureheight{0.45\textwidth*5/7}
% 	\setlength\figurewidth{0.45\textwidth}
% 	\InputIfFileExists{1_precom_l2_abs.tikz}{}{\textbf{!! Missing graphics !!}}
     \includegraphics[width=7cm]{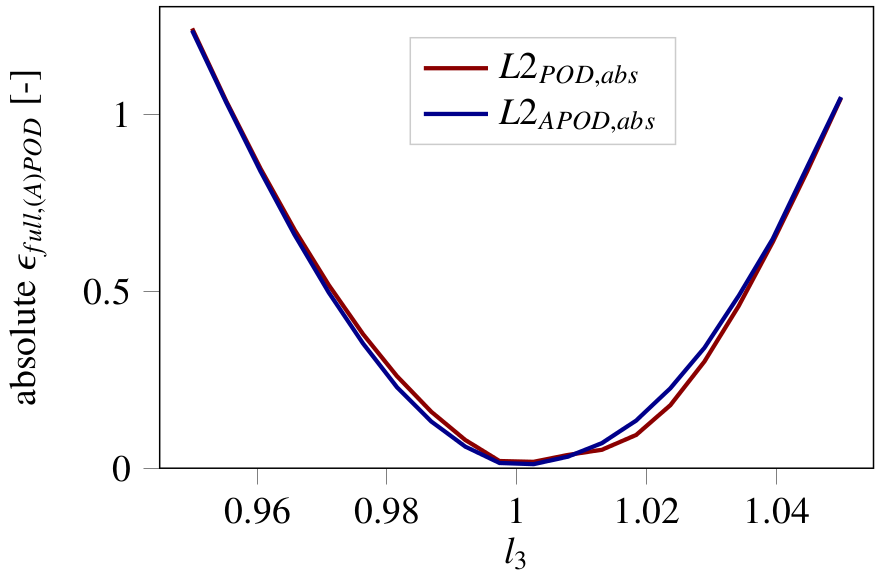}   
	\caption{ Comparison of the absolute $L2$ error }
	\label{fig:l2_abs}
%  \end{minipage}
%  \begin{minipage}{0.2\textwidth}
\end{figure}
\begin{figure}
	\centering
% 	\setlength\figureheight{0.45\textwidth*5/7}
% 	\setlength\figurewidth{0.45\textwidth}
% 	\InputIfFileExists{1_precom_qoi_abs.tikz}{}{\textbf{!! Missing graphics !!}}
     \includegraphics[width=7cm]{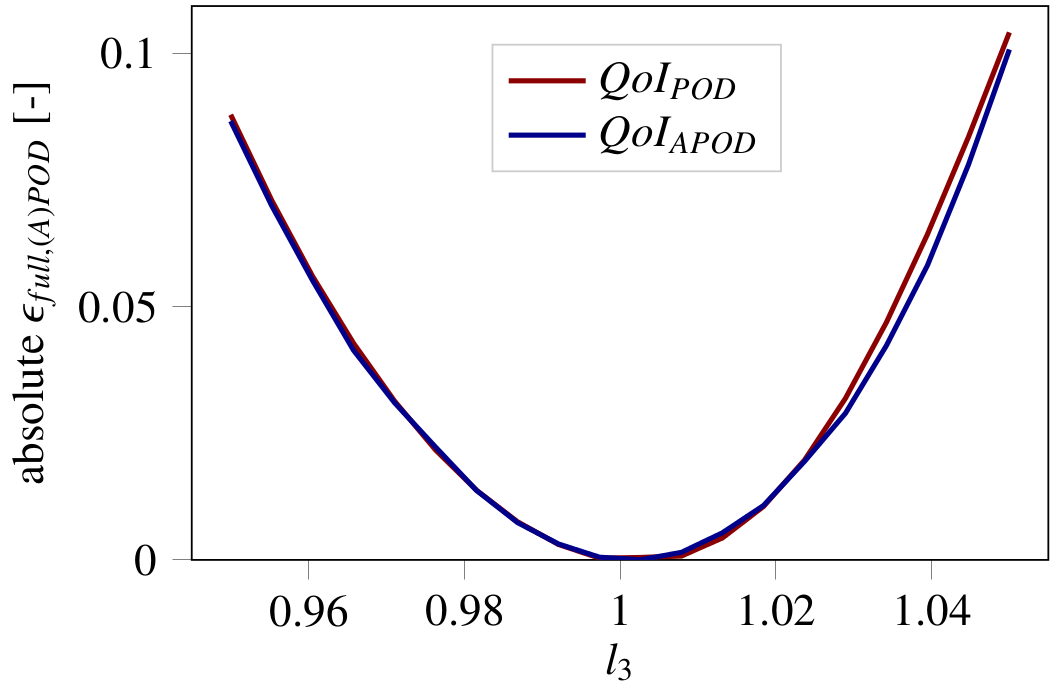}   
	\caption{ Comparison of the absolute QoI error }
	\label{fig:qoi_abs}
%  \end{minipage}
\end{figure}

POD and APOD yields accurate results if the Young's modulus is stiffer than the Young's modulus of the precalculation.
%Adding additional parameters leads to less accurate results. 
For our problem the dimensions of the cube are added as parameters. 
Each dimension varies between $l_1,l_2,l_3 \in [0.95,1.05]$. 
In Fig.~\ref{fig:l2_abs} and \ref{fig:qoi_abs} two points are observed. 
First, both POD and APOD are not able to capture the influence of change in the dimensions with the precomputation at hand. 
Second, both measures reflect this. 
Whereas, the measure via the QoI is easier to compute. 
Thus, in the following we will focus on the error measurement using the QoI.
To improve POD and APOD sensibly chosen precomputations are needed. 
Stiffer Young's moduli do not provide additional information, 
since these are equivalent to the effect of weaker loading and softer material,
which are already captured in the first time steps. 
For the dimensions the strongest effect on the quantity of interest is expected for 
the extremal values. 
Thus, the parameter configurations from Tab.~\ref{tab:cube_uncertain} are taken for the precalculations.

\begin{table}[ht]
 \caption{Material parameters and different sets of the precalculations}
 \label{tab:cube_uncertain}
 \begin{tabular}{lll}
  Material 	& Parameters 			& Pressure		\\
  Neo-Hookean	& $E=200 \SI{}{\N / \mm\squared}$		& $-0,320 \SI{}{\N / \mm\squared}$ 	\\
		& $\nu=0.4999$			&			\\ 
\hline
 Sets:		&				&			\\	
 \multicolumn{3}{l}{1.~~~~$l_1=0.95,~l_2=0.95,~l_3=0.95$} 	\\
 \multicolumn{3}{l}{2.~~~~$l_1=1.05,~l_2=0.95,~l_3=0.95$} 	\\
 \multicolumn{3}{l}{3.~~~~$l_1=0.95,~l_2=1.05,~l_3=0.95$} 	\\
 \multicolumn{3}{l}{4.~~~~$l_1=0.95,~l_2=0.95,~l_3=1.05$} 	\\
 \end{tabular}

\end{table}

\begin{figure}[ht]
% \begin{tikzpicture}
%     \begin{axis}[width=7cm, zmin=0, xlabel=$l_1$, ylabel=$l_3$, zlabel=$\epsilon_{full,(A)POD}$, colorbar, legend entries={POD, APOD}, legend image post style={scale=2}]
%     \addlegendimage{only marks, mark=square, color=red}
%     \addlegendimage{only marks, mark=square, color=blue}
%     \addplot3 [patch,patch refines=2,
% % 		shader=faceted interp,
% 		patch type=biquadratic,
% 		opacity=1] table{9_APOD.dat};
%     \addplot3[patch,patch refines=2,
% % 		shader=faceted interp,
% 		patch type=biquadratic,
% 		opacity=0.6] table{9_POD.dat};
%     \addplot3[patch,patch type=biquadratic,mesh,red, patch refines=2]
% 		 table{9_POD.dat};
%     \addplot3[patch,patch type=biquadratic,mesh,blue, patch refines=2]
% 		 table{9_APOD.dat};
% %     \addplot3 [only marks, mark options={scale=0.6}, red] file {extra_files/APOD/lange_095/E_200/coord_APOD.dat};
% %     \addplot3 [only marks, mark options={scale=0.6}, blue] file {extra_files/APOD/lange_095/E_200/coord_POD.dat};
%     \end{axis}
% \end{tikzpicture}
     \includegraphics[width=7cm]{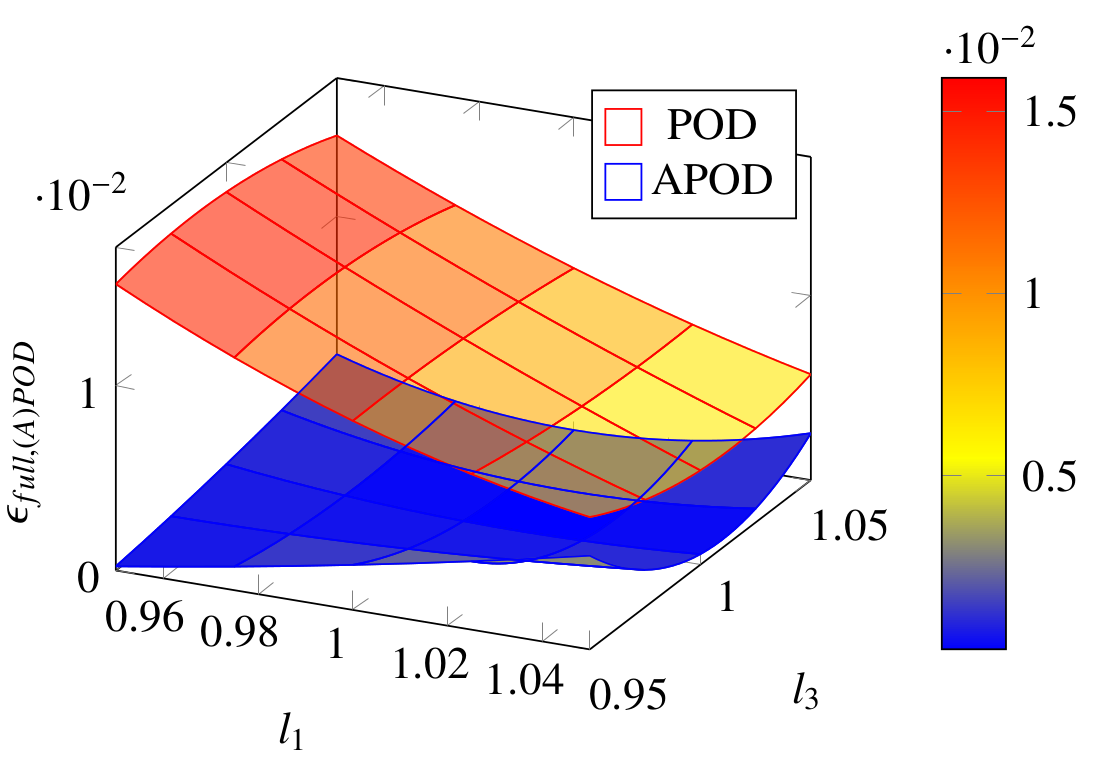}   
\caption{Comparison of the accuracy of POD and APOD for uncertain parameters: $l_2=1.00$, $E=200$ and $10$ modes}
\label{fig:3D_100_200}
\end{figure}

In Fig.~\ref{fig:3D_100_200} the results for the new precalculations are shown. 
The error $\epsilon_{\mathrm{full,(A)POD}}$ is plotted over the width $l_1$ and the height $l_3$ for $l_2=1.00$ and $E=200$. 
For all tested parameter configuration APOD has a higher accuracy.

$l_2=1$ is selected, since this choice is not covered by the precalculations directly. 
A variation of $l_2$ leads to similar results with the same set of snapshots.

% \begin{figure}[ht]
% \label{fig:3D_095_200}
% \begin{tikzpicture}
%     \begin{axis}[width=7cm, zmin=0, xlabel=$l_1$, ylabel=$l_3$, zlabel=Error, colorbar, legend entries={POD, APOD}, legend image post style={scale=2}]
%     \addlegendimage{only marks, mark=square, color=red}
%     \addlegendimage{only marks, mark=square, color=blue}
%     \addplot3 [patch,patch refines=2, % nodes near coord1={APOD},
% % 		shader=faceted interp,
% 		patch type=biquadratic,
% 		opacity=1] table{extra_files/APOD/lange_095/E_200/9_APOD.dat};
%     \addplot3[patch,patch refines=2,
% % 		shader=faceted interp,
% 		patch type=biquadratic,
% 		opacity=0.6] table{extra_files/APOD/lange_095/E_200/9_POD.dat};
%     \addplot3[patch,patch type=biquadratic,mesh,red, patch refines=2]
% 		 table{extra_files/APOD/lange_095/E_200/9_POD.dat};
%     \addplot3[patch,patch type=biquadratic,mesh,blue, patch refines=2]
% 		 table{extra_files/APOD/lange_095/E_200/9_APOD.dat};
% %     \addplot3 [only marks, mark options={scale=0.6}, red] file {extra_files/APOD/lange_095/E_200/coord_APOD.dat};
% %     \addplot3 [only marks, mark options={scale=0.6}, blue] file {extra_files/APOD/lange_095/E_200/coord_POD.dat};
%     \end{axis}
% \end{tikzpicture}
% \caption{Comparison of the accuracy of POD and APOD for uncertain parameters: $l_2=0.95$, $E=200$ and 10 modes}
% \end{figure}

% In order to see the influence of the parameters $l_2$ the results for a length $l_2=1.00$ and $l_3=1.05$ are plotted in Fig. \ref{fig:3D_100_200} and Fig. \ref{fig:3D_105_200}. The Young's modulus is with $E=200$ the same as in the precalculations. 

\begin{figure}[ht]
% \begin{tikzpicture}
%     \begin{axis}[width=7cm, zmin=0, xlabel=$l_1$, ylabel=$l_3$, zlabel=$\epsilon_{full,(A)POD}$, colorbar, legend entries={POD, APOD}, legend image post style={scale=2}]
%     \addlegendimage{only marks, mark=square, color=red}
%     \addlegendimage{only marks, mark=square, color=blue}
%     \addplot3 [patch,patch refines=2,
% % 		shader=faceted interp,
% 		patch type=biquadratic,
% 		opacity=1] table{9_APOD2.dat};
%     \addplot3[patch,patch refines=2,
% % 		shader=faceted interp,
% 		patch type=biquadratic,
% 		opacity=0.6] table{9_POD2.dat};
%     \addplot3[patch,patch type=biquadratic,mesh,red, patch refines=2]
% 		 table{9_POD2.dat};
%     \addplot3[patch,patch type=biquadratic,mesh,blue, patch refines=2]
% 		 table{9_APOD2.dat};
% %     \addplot3 [only marks, mark options={scale=0.6}, red] file {extra_files/APOD/lange_095/E_200/coord_APOD.dat};
% %     \addplot3 [only marks, mark options={scale=0.6}, blue] file {extra_files/APOD/lange_095/E_200/coord_POD.dat};
%     \end{axis}
% \end{tikzpicture}
     \includegraphics[width=7cm]{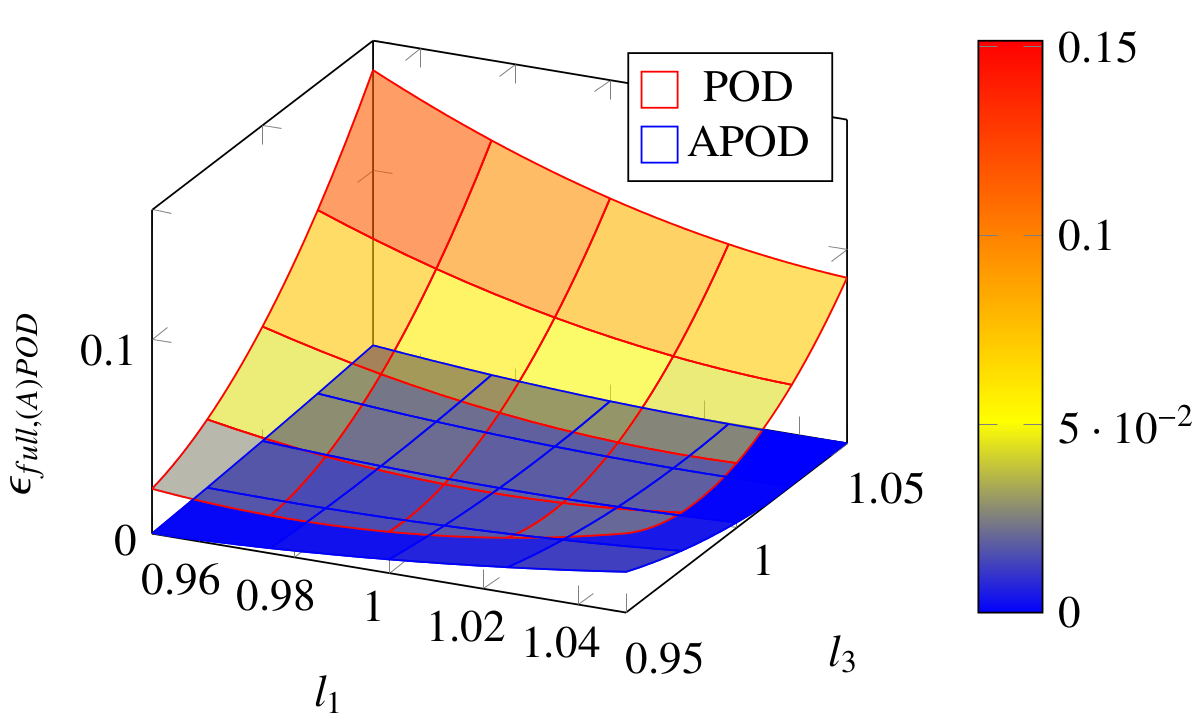}   
\caption{Comparison of the accuracy of POD and APOD for uncertain parameters: $l_2=1.05$, $E=2000$ and 10 modes}
\label{fig:3D_105_2000}
\end{figure}

\begin{figure}
% 	\centering
% 	\setlength\figureheight{0.45\textwidth*5/7}
% 	\setlength\figurewidth{0.45\textwidth}
% 	\InputIfFileExists{E-error.tikz}{}{\textbf{!! Missing graphics !!}}
     \includegraphics[width=7cm]{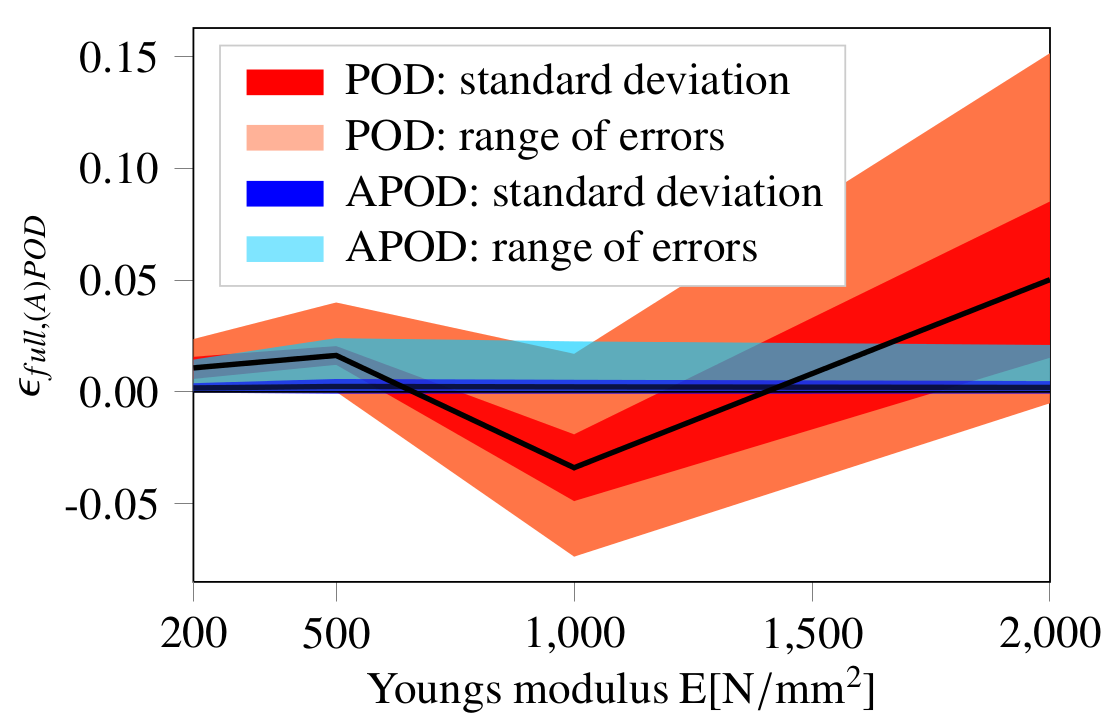}   
	\caption{Comparison of mean, deviation and range of errors for the range of uncertain parameters between POD and APOD}
	\label{fig:error_E_all}
\end{figure}

Raising Young's modulus to a value of $E=2000 \left[ \frac{N}{mm^2}\right]$ leads to worse accuracy of POD. This is observed in Fig.~\ref{fig:3D_105_2000}. 
The error $\epsilon_{\mathrm{full,POD}}$ rises up to $15\%$, while the maximal error $\epsilon_{\mathrm{full,APOD}}$ stagnates at $2\%$. 
It is therefore interesting to investigate the accuracy of POD and APOD for a variation of Young's modulus.
In Fig.~\ref{fig:error_E_all} this relation is presented. 
The modulus is fixed and the dimensions are covered by $1000$ samples.
The ranges and the standard deviation of the errors $\epsilon_{\mathrm{full,POD}}$ and $\epsilon_{\mathrm{full,APOD}}$ are plotted.
Both, range and standard deviation of $\epsilon_{\mathrm{full,APOD}}$ remain constantly smaller than the respective range and standard deviation of $\epsilon_{\mathrm{full,POD}}$.

% 
% \begin{figure}[ht]
% \caption{Comparison of the accuracy of POD and APOD for uncertain parameters: $l_2=1.05$, $E=1000$ and 10 modes}
% \label{fig:3D_105_1000}
% \begin{tikzpicture}
%     \begin{axis}[width=7cm, zmin=0, xlabel=$l_1$, ylabel=$l_3$, zlabel=Error, colorbar]
%     \addplot3 [patch,patch refines=2,
% % 		shader=faceted interp,
% 		patch type=biquadratic,
% 		opacity=1] table{extra_files/APOD/lange_105/E_1000/9_APOD.dat};
%     \addplot3[patch,patch refines=2,
% % 		shader=faceted interp,
% 		patch type=biquadratic,
% 		opacity=0.6] table{extra_files/APOD/lange_105/E_1000/9_POD.dat};
%    \addlegendentry{POD}
%    \addlegendentry{APOD}
% %     \addplot3 [only marks, mark options={scale=0.6}, red] file {extra_files/APOD/lange_095/E_200/coord_APOD.dat};
% %     \addplot3 [only marks, mark options={scale=0.6}, blue] file {extra_files/APOD/lange_095/E_200/coord_POD.dat};
%     \end{axis}
% \end{tikzpicture}
% \end{figure}

% error of 1000 samples
\begin{figure}
% 	\centering
% 	\setlength\figureheight{0.45\textwidth*5/7}
% 	\setlength\figurewidth{0.45\textwidth}
% 	\InputIfFileExists{error_histogram_pod_apod.tikz}{}{\textbf{!! Missing graphics !!}}
     \includegraphics[width=7cm]{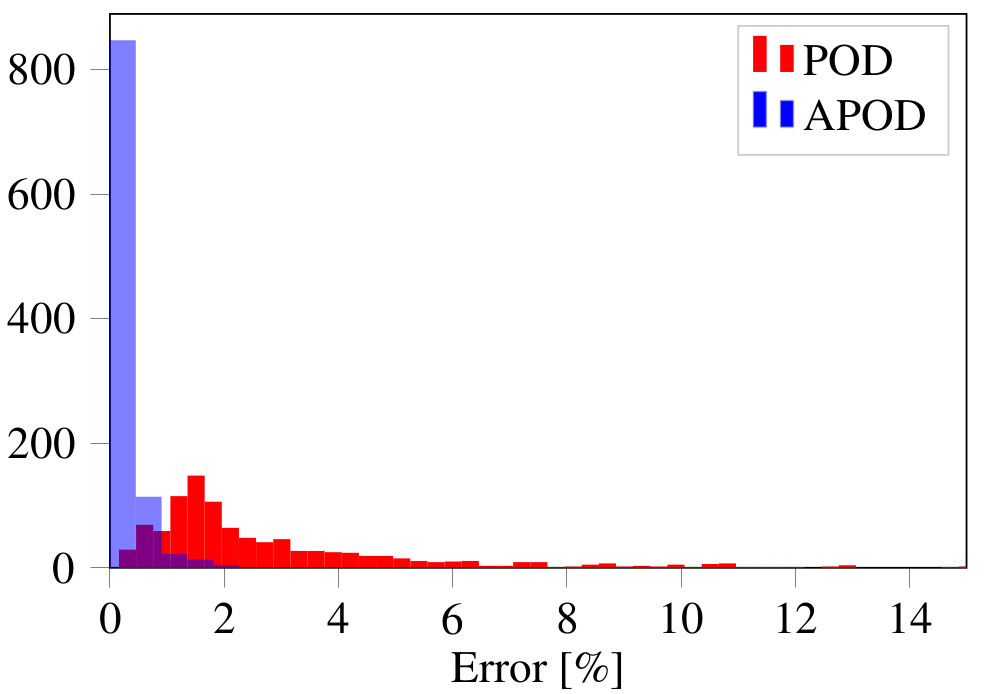}   
	\caption{Comparison of $\epsilon_{\mathrm{full,(A)POD}}$ for 1000 random samples.}
	\label{fig:error_POD_APOD_histo_1000}
\end{figure}

% Error for all calculations

Comparing the classic POD with APOD shows that the APOD increases the accuracy of the calculation significantly with the same or a smaller number of modes. 
This becomes particularly clear in Fig.~\ref{fig:error_POD_APOD_histo_1000}, where the error is measured for $1000$ random parameter configurations.
It requires a singular value decomposition for each time step. 
However, each SVD is cheaper than the SVD for all snapshots at once. 
By using APOD the computation cost of $\varphi$ is reduced drastically. 

Nonetheless, for accurate Monte Carlo simulations usually millions of samples are required,
thus the reduction of computational cost is not sufficient. 
In the next chapters a more specialized method and the synergies with (A)POD are investigated.

\section{Hierarchical tensor approximation (HTA)}
\label{sec:ht_theory}

One way to further reduce the complexity of the mapping from the high-dimensional parameter space to a manageable number of quantities of interest is the hierarchical Tucker decomposition \cite{grasedyck2010hierarchical,hackbusch2009new}.
To introduce this decomposition in a simple way, 
we discretise the parameter space $\mathcal{S}$ on a tensorized grid,
meaning that a finite number of points are taken for each parameter 
and every combination of these points are taken.
For our purposes $N$ equidistant points in each parameter direction are sufficient. 
Therefore, we denote for the parameters $l_i$ and $E$ the variations
\begin{align*}
	l_{i,k} = 0.95 + k \cdot \frac{0.1}{N-1} \; \text{and} \; E_k = 100 + k \cdot \frac{100}{N-1},
\end{align*}
with $k\in \left\{ 0,\ldots,N-1 \right\}$, and the resulting grid as $\mathcal{D}$.
We describe the corresponding indice set $\mathcal{I}$ by
\begin{align*}
	(i_1,i_2,i_3,i_4) \in 
	\mathcal{I} =
	\mathcal{I}_{ \left\{ 1 \right\} } \times \mathcal{I}_{ \left\{ 2 \right\} } \times 
	\mathcal{I}_{ \left\{ 3 \right\} } \times \mathcal{I}_{ \left\{ 4 \right\} }
	= \left\{ 0,1,\ldots,N-1 \right\}^{4}.
\end{align*}
Evaluating the mapping on this grid $\mathcal{D}$ yields a $4$-th order tensor array, 
where each entry $\Phi\left(i,j,k,l\right)$ is the maximal displacement 
for the parameter combination $(l_{1,i}, l_{2,j}, l_{3,k}, E_l)$.
Analogously, we denote the tensor arrays $\Phi_{\mathrm{POD}}$ and $\Phi_{\mathrm{APOD}}$ 
which stem from the prediction given by POD and APOD, respectively.

\begin{figure}
	\centering
	% Set the overall layout of the tree
% \tikzstyle{level 1}=[level distance=1.75cm, sibling distance=4cm]
% \tikzstyle{level 2}=[level distance=1.75cm, sibling distance=2.5cm]
% 
% % Define styles for bags and leafs
% \tikzstyle{root} = [text width=11em, text centered]
% \tikzstyle{inter} = [text width=12em, text centered]
% \tikzstyle{leaf} = [text width=9em, text centered]
% 
% \input{tree_modes}
% $ $ 
% \\[0.5cm]
% \input{tree_data}
     \includegraphics[width=6cm]{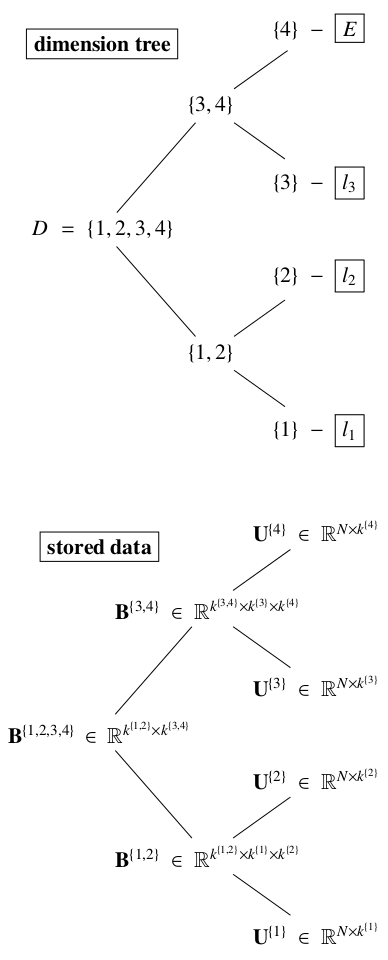}   
	\caption{The arrangement of the stored data in a hierarchical tucker format.}
	\label{fig:dimensiontree}
\end{figure}

For the HT decomposition of these tensor arrays each parameter is associated with a mode. 
For example the parameter $E$ is associated with the mode $\left\{ 4 \right\}$ 
and the parameters $l_1,l_2$ and $l_3$ with the modes $\left\{ 1 \right\},\left\{ 2 \right\}$ and $\left\{ 3 \right\}$.
We start with a set of modes $\left\{ 1,\ldots, d \right\}$. 
Then, this set of modes is split into $2$ subsets. 
The resulting subsets are recursively split until each subset only contains one mode.
Thus, the resulting structure is a binary tree with sets of modes in every vertex. 
This structure is denoted as \textit{dimension tree} $\mathcal{T}_D$, 
the root of this tree as $D=\left\{ 1,2,3,4 \right\}$ and the set of leaves as 
$ \mathcal{L} = \left\{ \left\{ 1 \right\}, \left\{ 2 \right\}, \left\{ 3 \right\}, \left\{ 4 \right\} \right\}$.
In Figure \ref{fig:dimensiontree} we find the dimension tree which is used throughout this paper.
A dimension tree is a concept to describe the structure of a HT tensor. 
When describing an algorithm for the HT format, 
the dimension tree is used to navigate us through the different parts of the HT tensor.

The second concept needed is the matricisation of a tensor, 
which is just a reordering of the tensor entries into a matrix.
The first two indices $(i_1,i_2)$  and the remaining two indices $(i_3,i_4)$ are renumbered 
lexicographically with a single index $i$ and a single index $j$,
 the matricisation of the tensor $A$ is denoted as
$$\mathcal{M}^{  \left\{ 1, 2 \right\}  }\left(A\right)_{i,j} \coloneqq A\left( i_1,i_2,i_3,i_4 \right). $$
For an accessible introduction to matricisations, see \cite{grasedyck2011introduction}.
In the same way, by assigning one single index $i$ to a set of indices corresoponding to the respective modes in the dimension tree 
and another indice $j$ to the remaining indices of the tensor, the matricisations is constructed for each node in the dimension tree. 
Therefore, the matricisation $\mathcal{M}^{ \left\{ 1,2,3,4 \right\} }(A)$ is assigned to the root. 
This matrix consists of a single column, which contains all tensor entries.
For the first leaf, the indice $i$ is associated with $(i_1)$ and the indice $j$ with $(i_2,i_3,i_4)$. 
This yields the matricisation $\mathcal{M}^{ \left\{ 1 \right\} }(A)$. 
In a way, this sets the first parameter in proportion to all remaining parameters. 
For example, if it holds 
$$ A(i_1,i_2,i_3,i_4) = c_{i_1},\quad \forall (i_2,i_3,i_4) \in \mathcal{I}, $$
the first parameter is unrelated to the rest. 
This leads to a rank $1$ matrix for the matricisation associated with this leaf.
The rank is an indication on how strongly one parameter relates to the remaining parameter.
With the association of the rank of a matricisation to the respective mode set $t \in \mathcal{T}_D$, 
it yields the definition of the \textit{hierarchical rank}
\begin{align*}
	\mathbf{k} = \left( k^t \right)_{t \in \mathcal{T}_D}, \quad k^t \coloneqq \mathrm{rank}\left( \mathcal{M}^t \left( A \right) \right), \quad t \in \mathcal{T}_D.
\end{align*}
The set of tensors with a certain hierarchical rank $\mathbf{k}$ is denoted as 
\begin{align*}
	\mathcal{H}_{\mathbf{k}} \coloneqq \left\{ A \in \mathbb{R}^{\mathcal{I}} : \mathrm{rank}\left( \mathcal{M}_t \left( A \right)  \right) \leq k_t, t\in \mathcal{T}_D \right\}.
\end{align*}

It is not to be expected that a tensor is included in $\mathcal{H}_{\mathbf{k}}$.
Therefore, the goal is rather, to find an approximation $\tilde{A} \in \mathcal{H}_{\mathbf{k}}$ for a tensor $A$.
In \cite{grasedyck2010hierarchical} an quasi-optimal algorithm was developed based on singular value decomposition (SVD).
For the resulting tensor $A_{\mathcal{H}-\mathrm{SVD}} \in \mathcal{H}_{\mathbf{k}}$ applies
\begin{align*}
	\| A - A_{\mathcal{H}-\mathrm{SVD}} \|_2 \leq \sqrt{2d-3} \min_{ A_{\mathrm{best}} \in \mathcal{H}_{\mathbf{k}} } \| A - A_{\mathrm{best}} \|_2 .
\end{align*}
The use of a SVD is expensive, in the sense that all entries of a matricisation are used. 
A data sparse method to find a low rank approximation of a matrix is the skeleton or cross approximation \cite{goreinov1997theory}.
The basic idea is to only use a subset of rows and columns of the matrix. 
Let $M\in \mathbb{R}^{\mathcal{I}_1 \times \mathcal{I}_2}$, $\mathcal{P} \subset \mathcal{I}_1$ and $\mathcal{Q} \subset \mathcal{I}_2$, then the cross approximation reads
\begin{align*}
	\tilde{M} \coloneqq M|_{\mathcal{I}_1 \times \mathcal{Q}} ~ S^{-1} ~ M|_{\mathcal{P} \times \mathcal{I}_2},
\end{align*}
where $S$ is the submatrix that arises at the intersections of the chosen rows and columns.
This method only employs the entries of the rows and columns, 
thus the computation of most of the entries in the matrix are spared.
Finding those pivot sets $\mathcal{P}$ and $\mathcal{Q}$ is possible 
in an adaptive fashion, cf.~\cite{bebendorf2000approximation}, meaning
that the algorithm proceeds until a cross approximation with a certain accuracy is found
or until a certain rank is reached.

With help of this algorithm, a low rank approximation is found
for each matricisation in the dimension tree. 
This yields column spaces for each node in $T_D$.
The important columns $ \left( U^{ \left\{ \mu \right\}}_{j} \right)_{j\in\left\{ 1, \ldots, k^{\mu} \right\}} $ 
found for the matricisations on the leafes are simply stored colum-wise in a matrix $ \mathbf{U}^{ \left\{ \mu \right\} } $, 
with $\mu \in \left\{ 1,2,3,4 \right\}$. We denote these as \textit{frames}.
One level higher, the important columns are found by the cross approximation method for a submatrix of $\mathcal{M}^{t}(A)$.
But instead of storing these columns, we compute and store the coefficients $\left(\mathbf{B}^{t}\right)_{i,j,k}$,
which approximate $U^{t}_{j}$ with linear combinations of kronecker 
products of the columns  $U^{s_1}_{k}\otimes U^{s_2}_{l}$,
with $\mathrm{sons}\left( t \right)=\left\{ s_1, s_2 \right\}$.
More precisely, these linear combinations read
\begin{align}
	\left( U^t \right)_{\cdot,j} = 
	\sum_{j_1  = 1}^{k_{s_{1}}} \sum_{j_2  = 1}^{k_{s_{2}}} 
	\left( \mathbf{B}^t \right)_{j, j_1, j_2} \left( U^{s_1} \right)_{\cdot, j_1} \otimes \left( U^{s_2} \right)_{\cdot, j_2}, \quad j \in \mathcal{I}^{t},
	\label{eq:recursive_ht}
\end{align}
for all $t\in \mathcal{T}_D \setminus \mathcal{L}$. 
Each transfer tensor $\mathbf{B}$ is computed by a series of projections between linear spaces. 
In condensed form it reads
\begin{align*}
\left( \mathbf{B}^{t} \right)_{j,\cdot,\cdot} = S^{-1}_{s_1} \left( \mathcal{M}^{t}(A) \right)_{(\cdot,\cdot),j} S^{-1}_{s_1}.
\end{align*}
With these \textit{transfer tensors} $\mathbf{B}^{t}$, the frames and Eq.~\eqref{eq:recursive_ht}, 
it yields a low rank approximation for a tensor array,
given that the ranks in each node are small. 
We denote this approximation in the following as hierarchical Tucker approximation (\textit{HTA})
and the algorithm for finding the HTA as \textit{generalized cross approximation}.
We have some control over the construction by setting the size of the submatrices and
by setting the accuracy or the maximal rank for each cross approximation.
This method of constructing a HT tensor is elaborated in \cite{ballani2013black}.

The transfer tensors and frames are the only data that is needed to be stored for an hierarchical tucker approximation.
Let $$k = \max_{t\in \mathcal{T}_{D}} k^{t},$$ then the upper bound reads $4 \cdot N \cdot k + 2k^3+k^2$
for the number of stored elements from Figure~\ref{fig:dimensiontree}. 
For $N=100$ and $k=5$, it yields $2275$ elements compared to $10^8$ elements in a full tensor.

With Eq.~\eqref{eq:recursive_ht} we have access to each entry of $U^{D}$
and therefore to each entry of the $\mathcal{H}$ tensor. 
For $(i_1,i_2,i_3,i_4) \in \mathcal{I}$ we pick the respective row $ v^{ \left\{ k \right\} } \coloneqq \mathbf{U}_{i_k,\cdot} $
and compute 
\begin{align*}
	\left( w^{ \left\{ 1,2 \right\} } \right)_{j} =  
	\sum_{j_1  = 1}^{k^{ \left\{ 1 \right\} }  } \sum_{j_2  = 1}^{k^{ \left\{ 2 \right\} }}
	\left( \mathbf{B}^{ \left\{ 1,2 \right\} } \right)_{j, j_1, j_2} 
	v^{ \left\{ 1 \right\} }_{j_1} v^{ \left\{ 2 \right\} }_{j_2}.
\end{align*}
In the same fashion, we compute $w^{ \left\{ 3,4 \right\} }$. 
Finally, the entry  is computed by
\begin{align*}
	\mathcal{H}( i_1,i_2,i_3,i_4 ) = 
	\sum_{j_1  = 1}^{k^{ \left\{ 1,2 \right\} }  } \sum_{j_2  = 1}^{k^{ \left\{ 3,4 \right\} }}
	\left( \mathbf{B}^{ \left\{ 1,2,3,4 \right\} } \right)_{j_1, j_2} 
	w^{ \left\{ 1,2 \right\} }_{j_1} w^{ \left\{ 3,4 \right\} }_{j_2}.
\end{align*}
Thus, this evaluation needs $2k+4$ matrix vector multiplications.

%\textcolor{red}{ Doing so, we implicitly introduce a form of regularization of the approximation.
%With lower ranks, we introduce higher regularity.}

% \setcounter{section}{5}
\section{HTA for the cube under compression }
\label{sec:HTA_cube}

%As briefly described in the previous chapter
%our model problem is the cube under compression, where 
%the maximal displacement in the middle of the cube is our quantity of interest.
Next, the mapping \eqref{eq:main_mapping} is approximated by the means of HTA.
For HTA an already available mapping is needed, e.g. $\varphi_{\mathrm{full}}$,
which is deduced from a full FEAP simulation.
Other choices are the reduced (A)POD models.
%Since we have access to a reduced model of the full simulation, we are also interested in the maximal displacement predicted by these reduced models 
%and the HTA which may be constructed from these reduced models.
The reduced mappings are denoted as $\varphi_{\mathrm{POD}}$ and $ \varphi_{\mathrm{APOD}} $ and
the respective approximating mappings in hierarchical Tucker format
are denoted as $\varphi_{\mathcal{H}_{\mathrm{full}}}, \varphi_{\mathcal{H}_{\mathrm{POD}}}$ and $ \varphi_{\mathcal{H}_{\mathrm{APOD}}} $.

We employ the generalized cross approximation as briefly described in Chapter \ref{sec:ht_theory} 
with random matrices of size $8$ for all tensor arrays, as long as the matricisation allows it.
If the number of rows or columns of the matricisation is smaller, the smaller number is used for the size of the respective random sub matrix. 
The cross approximation for each matricisation is then executed in an adaptive manner with a tolerance of $10^{-6}$
which means that the distance between the cross approximation and the actual matrix is less than $10^{-6}$ in each entry of the selected rows and columns.
This is not a guarantee for accuracy but at least a well established heuristic.

For each tensor array $\Phi$ which is deduced from the mapping $\varphi$, 
we examine different resolutions of the parameter space
and a different number of finite elements, namely $N \in \left\{ 4, 8, 16, 32, 64 \right\}$ 
and $N_{\mathrm{FEM}} \in \left\{ 2^3, 4^3, 8^3 \right\}$. 
To reduce computation time and to ensure the stability of the input for the HTA the loading is reduced by $80\%$.
The range of the Young's modulus is chosen to be in the range $\left[ 100,200 \right]$.

\begin{table}
	\centering
 \begin{tabular}{c|c|c|c|c|c|c}
    & \multicolumn{3}{c|}{Construction time [s]} & \multicolumn{3}{c}{Used entries}\\
    \hline
    sampling & \multicolumn{3}{c|}{mesh} & \multicolumn{3}{c}{mesh} \\
    rate     & $2^3$ & $4^3$ & $8^3$ & $2^3$ & $4^3$ & $8^3$ \\
    \hline
    4	& 44 & 272 & 5866 & 167 & 161 & 169 \\
    8   & 84 & 552 & 11135 & 338 & 341 & 335 \\
   16   & 123 & 798 & 16266 & 496 & 488 & 496 \\
   32   & 186 & 1230 & 24673 & 750 & 752 & 752 \\
   64   & 313 & 2066 & 41549 & 1264 & 1264 & 1264
%     4	& 44 & 272 & 5866 & 65.23 & 62.89 & 66.02 \\
%     8   & 84 & 552 & 11135 & 8.252 & 8.325 & 8.179\\
%    16   & 123 & 798 & 16266 & 0.757 & 0.745 & 0.757 \\
%    32   & 186 & 1230 & 24673 & 0.072 & 0.072 & 0.072 \\
%    64   & 313 & 2066 & 41549 & 0.008 & 0.008 & 0.008
 \end{tabular}
 \caption{Construction times of a HTA and the respective number of used entries.}
 \label{tab:build_hta}
\end{table}

% \begin{figure}
% 	\centering
% 	\setlength\figureheight{0.45\textwidth*5/7}
% 	\setlength\figurewidth{0.45\textwidth}
% 	\InputIfFileExists{extra_files/Plot/constr_time_vs_sampling.tikz}{}{\textbf{!! Missing graphics !!}}
% 	\caption{Construction times HTA.}
% 	\label{fig:constr_time}
% \end{figure}
% 
% \begin{figure}
% 	\centering
% 	\setlength\figureheight{0.45\textwidth*5/7}
% 	\setlength\figurewidth{0.45\textwidth}
% 	\InputIfFileExists{extra_files/Plot/used_entries_vs_sampling.tikz}{}{\textbf{!! Missing graphics !!}}
% % 	\includegraphics[width=7cm]{pictures/hta/used_entries_vs_sampling.png}
% 	\caption{Percentage of the needed entries needed to construct the HTA.}
% 	\label{fig:computed_entries}
% \end{figure}

The construction times are found for different numbers of samplings $N$ 
and $3$ different resolutions of the mesh are listed in Tab.~\ref{tab:build_hta}.
All computations were performed on a cluster with 2x Intel Xeon X5690 processors with 2x6 CPU cores with 3.47 GHz and 48GB RAM.
Each column shows the same trend. 
The number of finite elements determines the magnitude of the computational time needed, where
the sampling rate only changes the computational time in a sub-linear fashion. 
This is not a high price to pay, especially when compared to the number of used entries. 
Obviously, independently of the number of finite elements, 
the percentage of used entries of the tensor array drops rapidly with increasing sampling rate.
This means that a higher resolution in the parameter space leads to a higher compression 
of the encapsulated information by the HTA. 
Around $10^{3}$ evaluations are needed for $N=64$, independently of the number of finite elements used.
This indicates that the structure of the mapping $\phi$ is suitable for an approximation in a hierarchical Tucker format, 
stemming from simulations with different accuracies.

\begin{figure}
	\centering
% 	\setlength\figureheight{0.45\textwidth*5/7}
% 	\setlength\figurewidth{0.45\textwidth}
% % 	\InputIfFileExists{extra_files/Plot/min_max_error_2.tikz}{}{\textbf{!! Missing graphics !!}}
% 	\InputIfFileExists{min_max_error_stand_dev.tikz}{}{\textbf{!! Missing graphics !!}}
% % 	\includegraphics[width=7cm]{pictures/hta/err_minmax_el_8.png}
     \includegraphics[width=7cm]{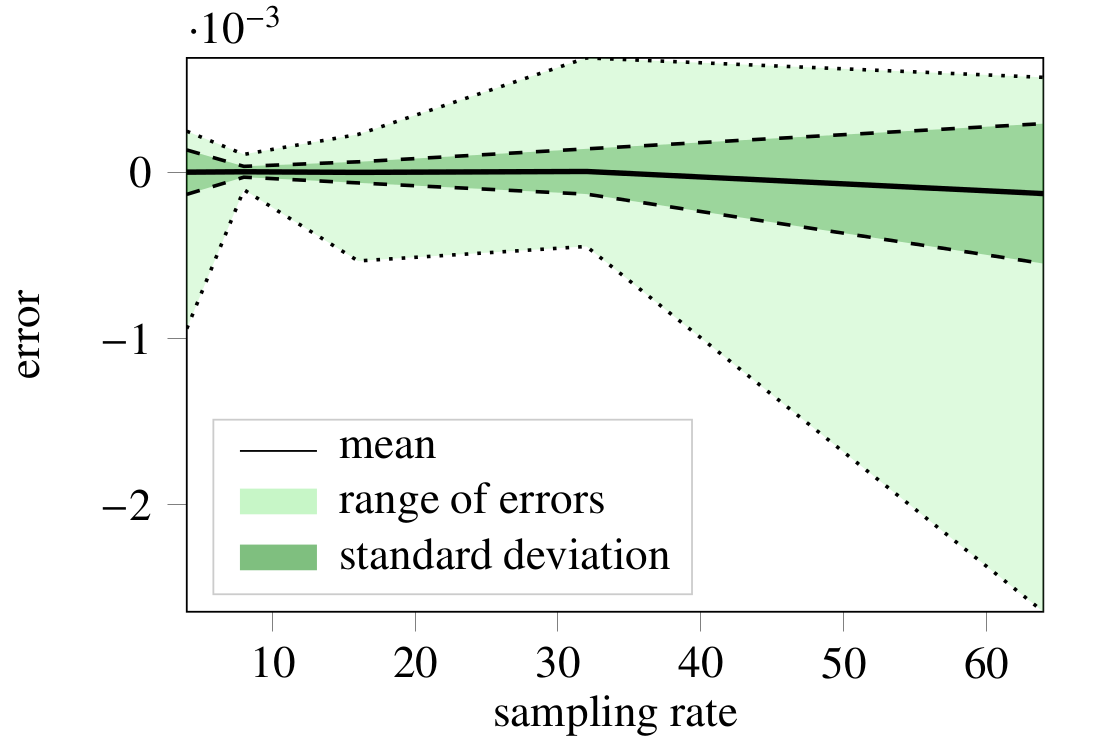}   
	\caption{Estimated absolute HTA error 
	$  \delta_{\mathrm{full},\mathcal{H}_{\mathrm{full}}} $ % = \left| \varphi_{\mathrm{F}}(l_1, l_2, l_3, E) - \mathcal{H}(l_1, l_2, l_3, E) \right| $  
	from $100$ random parameters
and $N_{\mathrm{FEM}} = 8^3 $.}
%\textcolor{red}{Warum nicht 1000?, Warum nicht zusätzlich mit std deviation wie bei den anderen Graphen dieser art?}}
	\label{fig:hta_err_minmax}
\end{figure}

\subsection{Approximation quality}
\label{sec:approx_quality_hta}
\begin{figure}
	\centering
% 	\setlength\figureheight{0.45\textwidth*5/7}
% 	\setlength\figurewidth{0.45\textwidth}
% %	\InputIfFileExists{extra_files/HTA_A-POD/HTA_unr_vs_unr.tikz}{}{\textbf{!! Missing graphics !!}}
% 	\InputIfFileExists{HTA_unr_vs_unr_1000.tikz}{}{\textbf{!! Missing graphics !!}}
     \includegraphics[width=7cm]{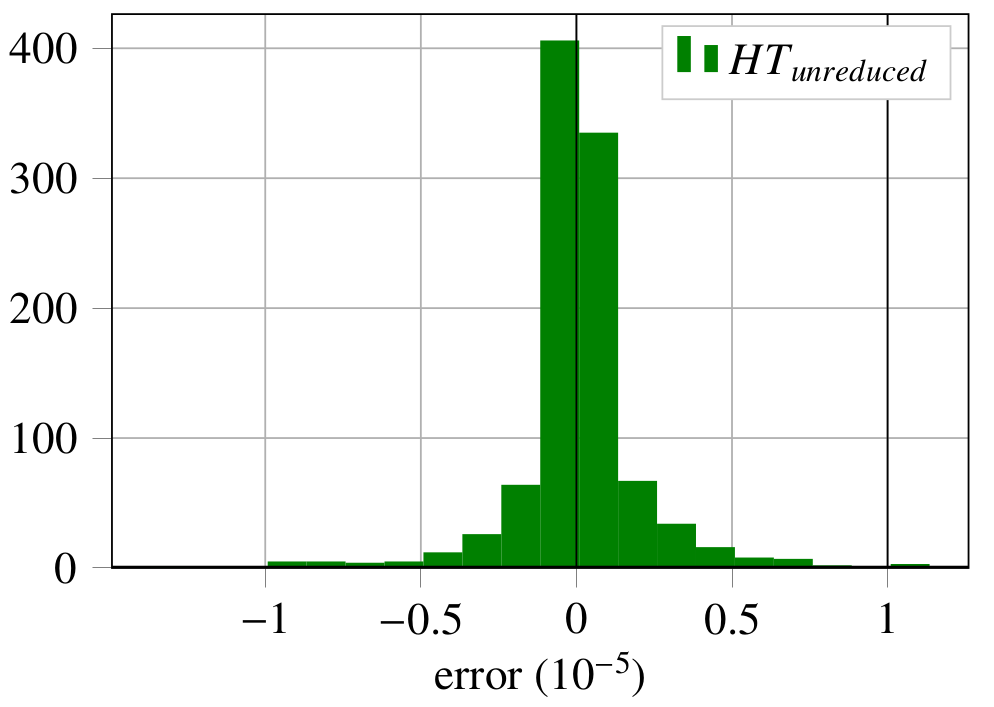}
	\caption{
		Distribution of all relative errors $\epsilon_{\mathrm{full}, \mathcal{H}_{\mathrm{full}}}$.
		Estimated from $1000$ random parameters on the tensorized grid.}
	\label{fig:error_hta_unr_unr}
\end{figure}

\begin{figure}
	\centering
% 	\setlength\figureheight{0.45\textwidth*5/7}
% 	\setlength\figurewidth{0.45\textwidth}
% 	\InputIfFileExists{eps_full_a-pod_easy.tikz}{}{\textbf{!! Missing graphics !!}}
     \includegraphics[width=7cm]{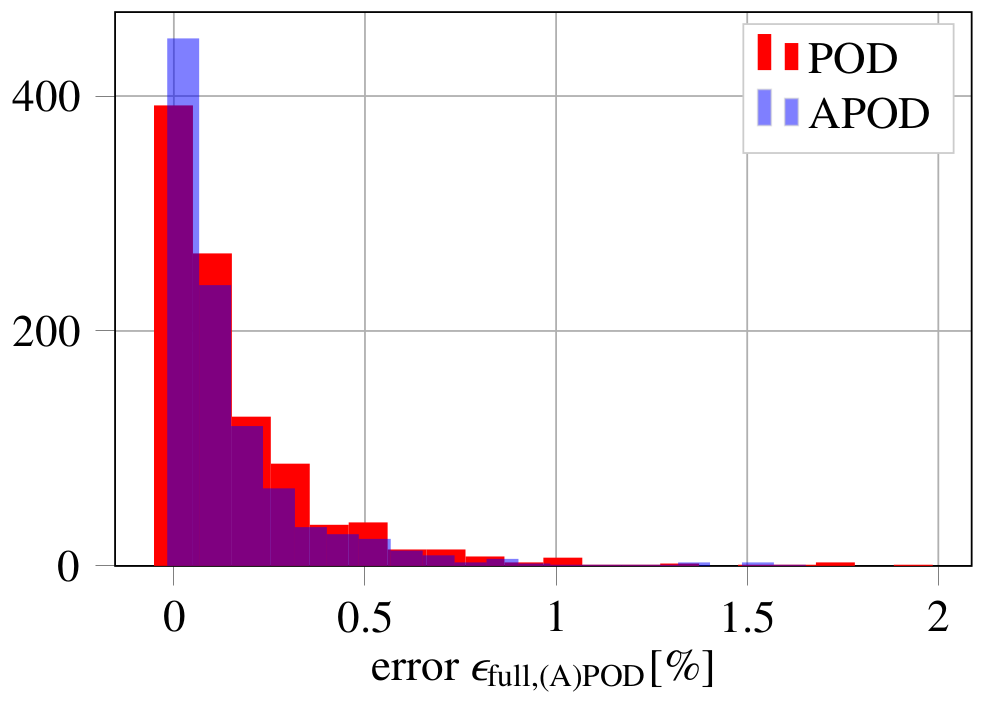}
	\caption{
		Distribution of all relative errors $\epsilon_{\mathrm{full}, \mathrm{(A)POD}}$.
		Estimated from $1000$ random parameters on the tensorized grid.}
	\label{fig:error_full_A-POD_easy}
\end{figure}

\begin{figure}
	\centering
% 	\setlength\figureheight{0.45\textwidth*5/7}
% 	\setlength\figurewidth{0.45\textwidth}
% % 	\InputIfFileExists{extra_files/HTA_A-POD/HTA_A-POD_vs_A-POD.tikz}{}{\textbf{!! Missing graphics !!}}
% 	\InputIfFileExists{HTA_A-POD_vs_A-POD_1000.tikz}{}{\textbf{!! Missing graphics !!}}
     \includegraphics[width=7cm]{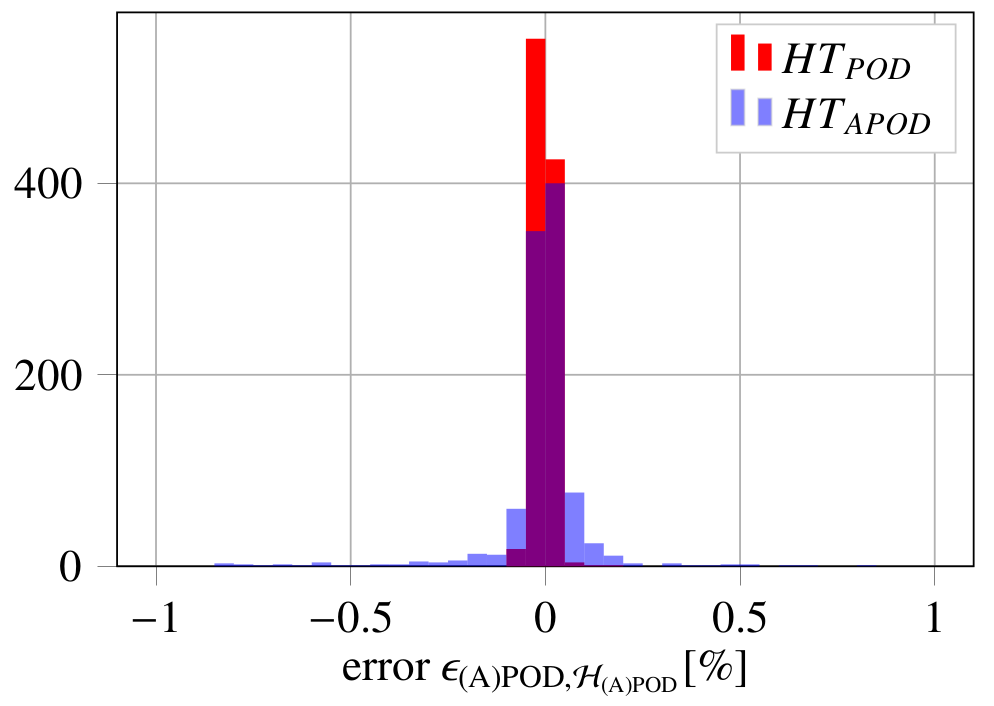}
	\caption{ Distribution of all relative errors $\epsilon_{\mathrm{APOD}, \mathcal{H}_{\mathrm{(A)POD}}}$.
		  Estimated from $1000$ random parameters on the tensorized grid.}
	\label{fig:error_hta_a-pod_a-pod}
\end{figure}

\begin{figure}
	\centering
% 	\setlength\figureheight{0.45\textwidth*5/7}
% 	\setlength\figurewidth{0.45\textwidth}
% % 	\InputIfFileExists{extra_files/HTA_A-POD/HTA_A-POD_vs_unr.tikz}{}{\textbf{!! Missing graphics !!}}
% 	\InputIfFileExists{HTA_A-POD_vs_unr_1000.tikz}{}{\textbf{!! Missing graphics !!}}
     \includegraphics[width=7cm]{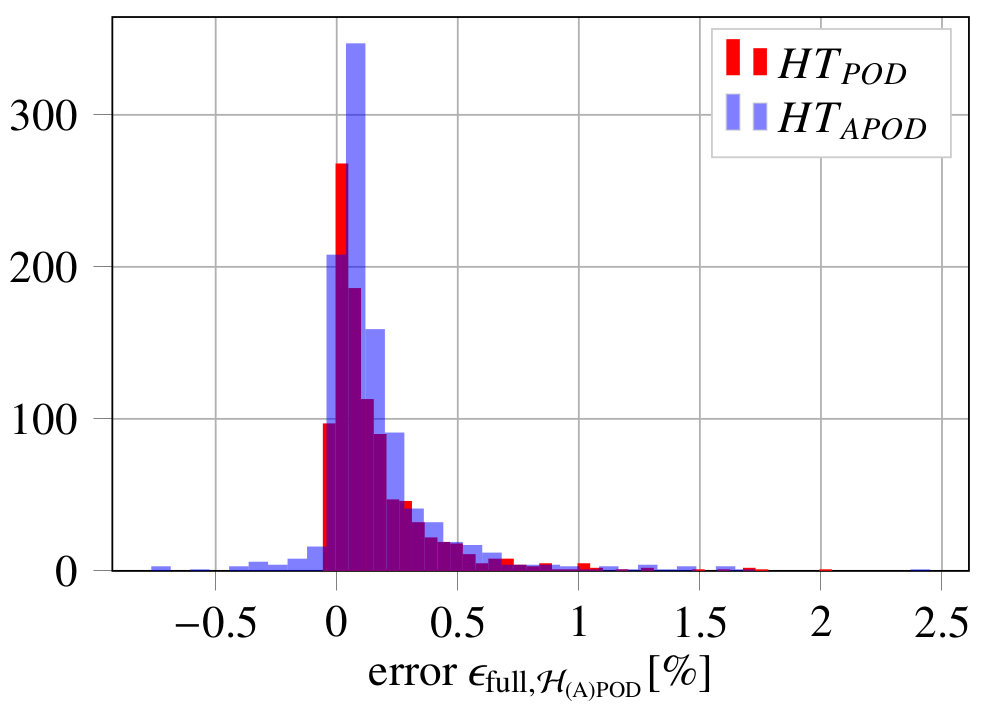}
	\caption{ Distribution of all relative errors $\epsilon_{\mathrm{full}, \mathcal{H}_{\mathrm{APOD}}}$ .
		  Estimated from $1000$ random parameters on the tensorized grid.
		}
	\label{fig:error_hta_a-pod_unr}
\end{figure}

\begin{table}
 \begin{tabular}{c|c|c|c|c}
  &       & std       &     & \\
  & mean  & deviation & max & min \\
  & [\%]  & [\%] & [\%] & [\%] \\
  \hline
  $\epsilon_{\mathrm{full,POD}}$ & 0.1641 & 0.2372 & 1.9854 & -0.0519 \\
  $\epsilon_{\mathrm{full,APOD}}$ & 0.1543 & 0.2047 & 1.6533 & -0.0173 \\
  $\epsilon_{\mathrm{POD}, \mathcal{H}_{\mathrm{POD}}}$ & 0.0030 & 0.0189 & 0.1586 & -0.0712 \\
  $\epsilon_{\mathrm{APOD}, \mathcal{H}_{\mathrm{APOD}}}$ & 0.0012 & 0.1618 & 2.4481 & -0.8481 \\
  $\epsilon_{\mathrm{full}, \mathcal{H}_{\mathrm{full}}}$ & 0.0000 & 0.0002 & 0.0011 & -0.0014 \\
  $\epsilon_{\mathrm{full}, \mathcal{H}_{\mathrm{POD}}}$ & 0.1611 & 0.2316 & 2.0459 & -0.0573 \\
  $\epsilon_{\mathrm{full}, \mathcal{H}_{\mathrm{APOD}}}$ & 0.1555 & 0.2506 & 2.4506 & -0.766 
 \end{tabular}
 \caption{Errors for different methods from $1000$ random parameters on the tensorized grid.}
 \label{tab:errors_all}
\end{table}

More interesting is the approximation quality of the HTA.
To estimate the quality we examine the distributions of various errors. 
For each distribution $1000$ random parameter combinations $(l_1,l_2,l_3,E) \in \mathcal{D}$ are used.
% First the absolute errors are observed.
%$$  \delta_{\mathrm{full},\mathcal{H}_{\mathrm{full}}} = \left| \varphi_{\mathrm{F}}(l_1, l_2, l_3, E) - \varphi_{\mathcal{H}(l_1, l_2, l_3, E)} \right|.$$ 
In Fig.~\ref{fig:hta_err_minmax} the range of errors is plotted including the mean over the different sampling rates and $N_{\mathrm{FEM}}=8^{3}$.
It is observed that the maximal error lies around the magnitude $10^{-4}$. 
In relation to the displacement, which lies between $20\%$ and $30\%$, the maximal error is six magnitudes smaller. 
%We also see that the sampling rate has only a small influence on the quality of the HTA for the unreduced model.
Accordingly, it is observed in Fig.~\ref{fig:error_hta_unr_unr}  that the maximal relative error 
$$\epsilon_{\mathrm{full},\mathcal{H}_{\mathrm{full}}}(p) = \frac{\varphi_{\mathrm{full}}(p) - \varphi_{\mathcal{H}_{\mathrm{full}}}(p) }{\varphi_{\mathrm{full}}(p)}$$
is of magnitude $10^{-6}$ in most cases.

%\textcolor{red}{Wie ist das hier? Haben wir erklärt wo die fehler herkommen nehmen wir nur eine Feap auflösung und sampling rate?}

The next error observed is the relative error $$\epsilon_{\mathrm{(A)POD}}(p) = \frac{\varphi_{\mathrm{(A)POD}}(p) - \mathcal{H}_{\mathrm{F}}(p) }{\varphi_{\mathrm{(A)POD}}(p)}.$$
Here,  Fig.~\ref{fig:error_hta_a-pod_a-pod} shows that the maximal relative error is smaller than $0.15\%$, which is a satisfactory result for certain applications. 
%\textcolor{red}{klären:
It suggests that the reduced model produces a mapping, which is less suitable for this kind of low rank approximation.
In comparison to the HTA of the unreduced model the error is $3$ orders of magnitude worse. 
%}
%\textcolor{red}{Brauchen wir hier an der Stelle nun nicht ein histogram mit dem relativen fehler von  $\epsilon_{\mathrm{(A)POD}}(p) = \frac{\varphi_{\mathrm{F}}(p) - \varphi_{\mathrm{F}}(p) }{\varphi_{\mathrm{F}}(p)}$. }

Whenthe HTA of a reduced model to approximate the mapping $\varphi_{\mathrm{F}}$ is used, 
the approximation errors of the reduced model and of the HTA are added. 
Empirically, this is observed with the relative error 
$\epsilon_{\mathrm{full},\mathcal{H}_{ \mathrm{(A)POD }} }(p)$,% = \frac{\varphi_{ \mathrm{full} }(p) - \varphi_{\mathcal{H}_{ \mathrm{(A)POD }}}(p) }{\varphi_{\mathrm{full}}(p)}$,
which is found in Fig.~\ref{fig:error_hta_a-pod_unr}. 
A majority of relative errors are smaller than $1\%$,  which is still a sufficient approximation quality for some applications.

An overview of all relative errors is found in Table~\ref{tab:errors_all}. 
APOD yields the better results, although it scatters more. 
Note that these errors were measured for a fifth of the loading. 
Therefore, the errors $\epsilon_{\mathrm{full,POD}}$ and $\epsilon_{\mathrm{full,APOD}}$ 
do not comply with the errors found in Sec.~\ref{sec:APOD}.
With a fifth of the original loading the cube shows less pronounced non-linear behaviour.
Therefore, to obtain the same error level, about the same number of modes are used for POD and APOD.

The HTA is seemingly better suited for POD than for APOD.
If the HTA is constructed from (A)POD,  the mean relative error to the unreduced model is
just the sum of the means of $\epsilon_{\mathrm{POD}, \mathcal{H}_{\mathrm{(A)POD}}}$ and $\epsilon_{\mathrm{full}, \mathcal{H}_{\mathrm{(A)POD}}}$.
The relative errors of HTA and the (A)POD approximation are two magnitudes apart from each other.
Overall, it is feasible to surrogate $\varphi_{\mathrm{(A)POD}}$ with $\varphi_{\mathcal{H}_{\mathrm{(A)POD}}}$.

Of course, there is plenty of room for improvement and many angles to approach this insufficiency.
It would be easily possible to increase the ranks of the HTA or the number of modes and snapshots of the POD.
Both ideas would lead to an increase of computational cost,
but the balance between these cost and the approximation quality is its own task for every application, 
where HTA and POD methods are used in combination.

Instead of constructing the approximation  $\mathcal{H}_{\mathrm{(A)POD}}$ of the (A)POD model we will employ the HTA to improve the used snapshots and therefore the (A)POD model itself in the next chapter.
This is to our knowledge the first attempt in literature.
%We will attempt to improve the choice of snapshots by using a series of HTA approximations.

\subsection{Finding snapshots with HTA}
\label{sec:snapshot_search}
%
%	Das ist klug, dass wir den l2 Fehler berechnet haben 
%	Weil damit die Idee aufgeht für eine Größe (displacement) die residueen zu vergelichen
%	und damit aber für eine andere Größe (l2 fehler) auch die Approximationseigenschaften vergleichen.
%

In the previous chapters, the POD and APOD model was built upon snapshots, 
which were carefully chosen from an understanding of the model problem.
Therefore, one question at hand is, if it is possible to improve the POD or APOD method, 
by finding better snapshots with the help of the HTA of the unreduced method. 
%Especially, since 
Since the approximation quality of HTA for the unreduced model has a high accuracy,
the CPU time effort can be significantly reduced by using HTA instead of the full FE model.
%, when we compare the maximal displacement of the reduced and unreduced model. 

One way to assess the approximation quality of the (A)POD model is the comparison to the unreduced FE model.
With this measure at hand it is possible to find parameter configurations where the discrepancy between FE model and (A)POD model
is especially large. New snapshots are then generated for this particular parameter configurations, thus reducing the discrepancy between FE and (A)POD model.
One iteration of this procedure may be be written down as 
\begin{itemize}
	\item[(1)] Build a (A)POD model $\varphi_{\mathrm{(A)POD},k}(p)$ from an initial set of snapshots $\mathcal{S}_k$
	\item[(2)] Provide a measure $r: \mathcal{D} \to \mathbb{R}$ for the discrepancy between reduced and unreduced model
	\item[(3)] Find $p^* =\underset{p\in\mathcal{D}}{\operatorname{argmax~}} |r(p)| $
	\item[(4)] Generate the snapshots $\mathcal{S}(p^*)$ and enrich the used set of snapshot $\mathcal{S}_{k+1} = \mathcal{S}_k \cup \mathcal{S}(p^*)$
\end{itemize}
This step may be repeated until a fixed number of snapshots is found or a certain tolerance is reached.
Step $3$ is especially expensive, since the parameter space may be high dimensional. 
Additionally, for an accurate measure of the discrepancy evaluations of the full model are needed.
To mitigate these costs, we propose to use HTA as a surrogate model of the residual $r$ in each step.
Therefore, let the residual be denoted as
\begin{align}
	r_k: \mathcal{D} \to \mathbb{R}: p \mapsto \varphi_{\mathcal{H}_{\mathrm{full}}}(p) -  \varphi_{\mathrm{(A)POD},k}(p),
	\label{eq:residual}
\end{align}
where $\varphi_{\mathrm{(A)POD},k}(p)$ denotes the $k$-th reduced model from a successively constructed set of models.
Naturally, we start this series with $\varphi_{\mathrm{(A)POD},0}(p)\equiv0$ which is equivalent to an empty set of snapshots $\mathcal{S}_0 = \emptyset$ and implies $r_0(p)=\varphi_{\mathcal{H}_{\mathrm{full}}}(p)$.
The main step is the third step which consists of the computation of $\max_{p \in \mathcal{D}} \left | r_k(p) \right |$.

The maximum is estimated by finding the maxima in alternating directions iteratively.
Heuristically, in our case this gives a good estimation for the maximum.
To improve the snapshot set it is not necessary to find the global maximum.

% The maximum is estimated by the following procedure.
% Starting from an arbitrary $p \in \mathcal{D}$ all directions are fixed but the first, 
% then the absolute maximum in this direction is found. 
% At this new point all directions but the second  are fixed and again the absolute maximum found. 
% This steps are repeated iteratively until the found value does not change more than a tolerance of $10^{-6}$. 
% Heuristically, in our case this gives a good estimate for the absolute maximum.

Once a maximum (global or local) is determined, this leads to a parameter configuration where the discrepancy between reduced and unreduced model is high.
Exactly this parameter configuration is then used to generate a new set of snapshots $\mathcal{S}_{k+1}$, leading to $\varphi_{\mathrm{(A)POD},k+1}$.
These steps are repeated iteratively until a certain tolerance is fulfilled or a maximum number of iterations is reached.

% 
% The reduced model for $\varphi_{\mathrm{(A)POD},1}(p)$ is built from these snapshots.
% Next, a new  HTA is constructed for $r_1$. 
% Note that this step is not necessary, since, in our case, it is also feasible to use the HTA and the reduced model 
% to compute the absolute maximum of $r_1$. 
% But for more sophisticated models, where the evaluation of the reduced model is also sufficiently costly, 
% this step is adequate.
% Again, we estimate the absolute maximum and the corresponding parameter configuration, use this parameter configuration
% to generate additional snapshots, which are used to enrich the preceding reduced POD model.
% This procedure is repeated several times. 

We perform this procedure for POD and APOD, with underlying a FE model with $4^3$ elements 
and use a HTA with a sampling rate of $64$ points per parameter.
In our experiments this coarse FE model lead to very similar parameter configurations as a finer model.
This strategy makes the procedure of this section feasible in practise to provide a sensible selection of snapshots for more accurate models.

In Tab.~\ref{tab:r_k_progress}, the devolopment of the absolute maximum with each iteration is found.
The decline of the residual is significant until it stagnates after the fifth iteration. 
This, most probably, leads back to the approximation error of the HTA for the residual $r_k$.
%\textcolor{red}{
%Dann müsste man vermutlich das maximum von $r_k$ nochmal anders bestimmen.
%Also das maximum ohne HTA von rk suchen.
%}
\begin{table}
	\centering
	\begin{tabular}{c|c|c}
		k & POD & APOD  \\
		\hline
		0 & 34.628   & 34.6262 \\
		1 & 11.8444   & 1.52712e+08 \\
		2 & 0.145506   & 0.0933674 \\
		3 & 0.206719   & 0.178969 \\
		4 & 0.0496022  & 0.0337002 \\
		5 & 0.0352917 & 0.0618711 
% Number samples: 64, Elements: 2x2x2  -> High samples -> Low Calculation time		
	\end{tabular}
	\caption{The development of $\max_{p \in \mathcal{D}} \left| r_k(p) \right|$. }
	\label{tab:r_k_progress}
\end{table}

\begin{figure}
	\centering
% 	\setlength\figureheight{0.45\textwidth*5/7}
% 	\setlength\figurewidth{0.45\textwidth}
% 	\InputIfFileExists{l2-HTA.tikz}{}{\textbf{!! Missing graphics !!}}
     \includegraphics[width=7cm]{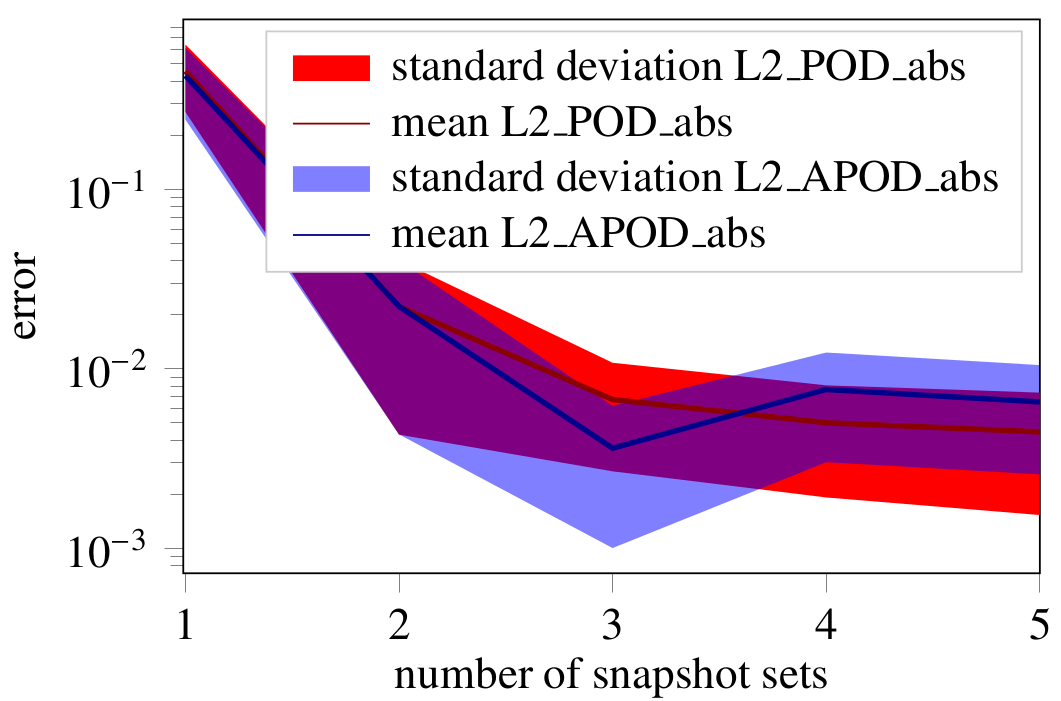}
	\caption{%Zeigt wie der l2 fehler mit mehr hta gewählten snapshots sich so verhält.
		Development of the $L_2$-error between the full model and the (A)POD model with a successively enriched snapshot set.
	}
	\label{fig:quality_hta_snapshots_l2}
\end{figure}

\begin{figure}
	\centering
% 	\setlength\figureheight{0.45\textwidth*5/7}
% 	\setlength\figurewidth{0.45\textwidth}
% 	\InputIfFileExists{qoi-HTA.tikz}{}{\textbf{!! Missing graphics !!}}
     \includegraphics[width=7cm]{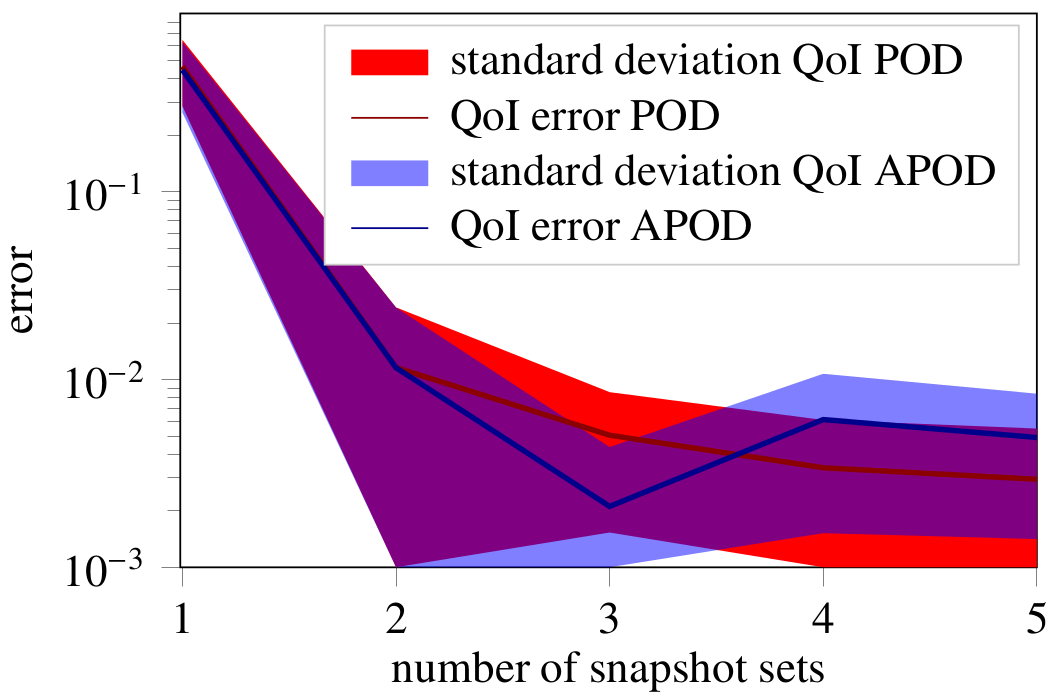}
	\caption{%Zeigt wie der qoi fehler mit mehr hta gewählten snapshots sich so verhält.
		Development of the QoI-error between the full model and the (A)POD model with a successively enriched snapshot set.
	}
	\label{fig:quality_hta_snapshots_qoi}
\end{figure}

\begin{figure}
	\centering
% 	\setlength\figureheight{0.45\textwidth*5/7}
% 	\setlength\figurewidth{0.45\textwidth}
% 	\InputIfFileExists{qoi-HTA-vs-MAN-APOD.tikz}{}{\textbf{!! Missing graphics !!}}
     \includegraphics[width=7cm]{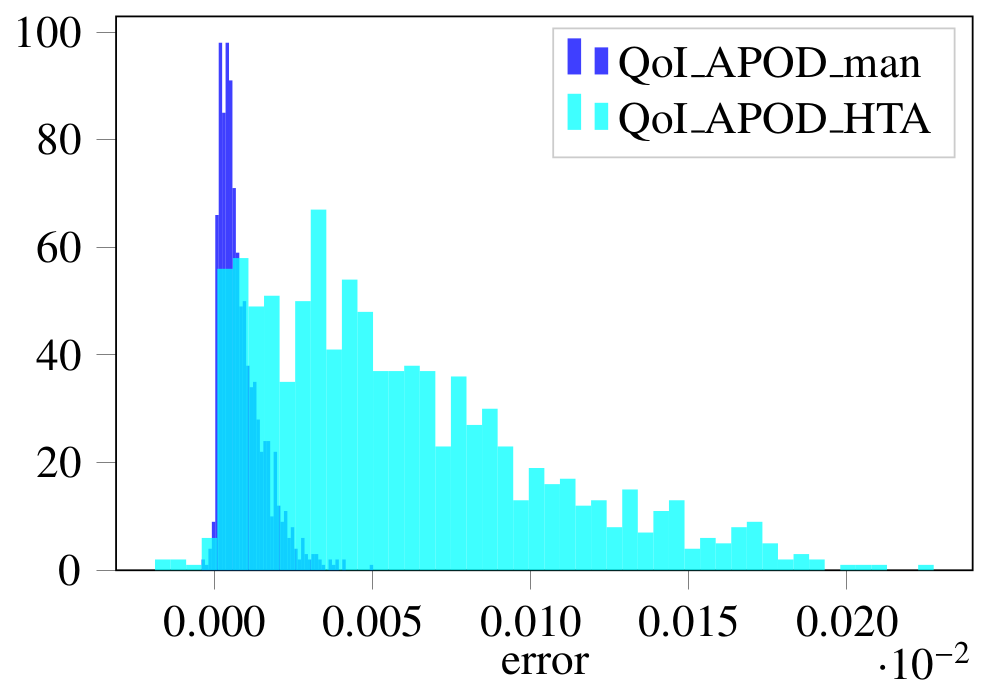}
	\caption{%Zeigt den qualitätsunterschied zwischen manuellen und hta gewählten snapshots.
	Distribution of the errors for two different APOD models from $1000$ samples. 
	With manually chosen snapshots and with systematically found snapshots using HTA.
	}
	\label{fig:comparison_man_vs_hta_snapshots}
\end{figure}

To compare these new models with the models in Sec.~\ref{sec:APOD}, 
we generate random parameters and compute the difference in the quantity of interest as one measurement
for the quality of the reduced model, as well as the $L^2$ error as a second measurement. 
We observe in Fig.~\ref{fig:quality_hta_snapshots_l2} and Fig.~\ref{fig:quality_hta_snapshots_qoi}, 
that the whole procedure is working for both measurements.
For APOD we find that there is a kink for both measurements.
On the other hand, a comparison with the manually chosen
snapshots in Fig.~\ref{fig:comparison_man_vs_hta_snapshots} shows, 
that the HTA chosen snapshots are significantly worse for APOD. 
One reason could be, that the interesting parameter configurations are found
at the boundary of the parameter space, but the HTA lacks accuracy at the boundaries.
The reason for this is the manner in which the HTA is constructed. 
Since the generalized cross approximation takes submatrices of the respective matricizations, 
the parameters on the boundary space are usually not hit upon.
For this simple problem it is practicable to choose suitable snapshots. 
In more challenging situations this task becomes more difficult. 
In those cases, the procedure above could suggest parameter configurations for a set of quantities of interest,
so that the choice of snapshots becomes feasible.

\begin{figure}
	\centering
% 	\setlength\figureheight{0.45\textwidth*5/7}
% 	\setlength\figurewidth{0.45\textwidth}
% 	\InputIfFileExists{qoi-HTA-vs-MAN-POD.tikz}{}{\textbf{!! Missing graphics !!}}
     \includegraphics[width=7cm]{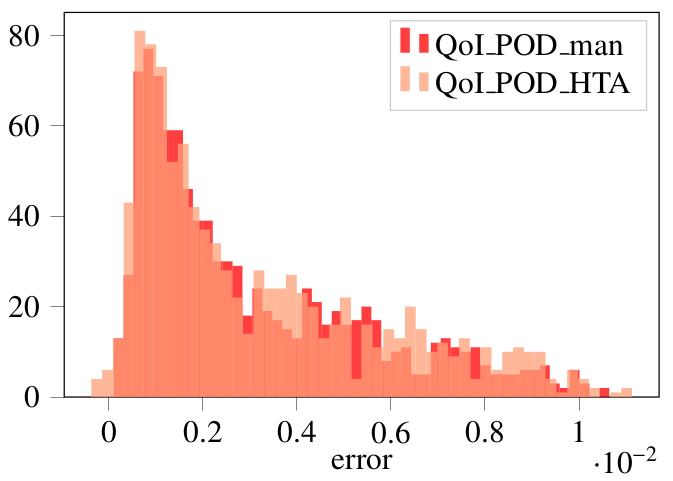}
	\caption{%Zeigt den qualitätsunterschied zwischen manuellen und hta gewählten snapshots.
		Distribution of the errors for two different POD models from $1000$ samples. 
		With manually chosen snapshots and with systematically found snapshots using HTA.
	}
	\label{fig:comparison_man_vs_hta_snapshots-POD}
\end{figure}

%\textcolor{red}{ Wie sah nochmal der Vergleich POD man gegen POD HTA aus? }

Once a suitable reduced model is found, it may be used for a variety of application.
Since the main application in this paper is uncertainty quantification in form of a Monte Carlo estimation of the mean and variance,
we will discover in the following chapter how to perform this with the HTA.

\section{ Comparison - HTA and APOD }
\label{sec:comparison_HTA_APOD}

In this paper, three levels of modeling are considered.
The first level is the FEM model, the second level is the POD model
deduced from the FEM model and the third level is the low rank approximation
in the hierarchical tensor format. 
Each level is in some sense more specialized and carries less information. 
E.g., the FEM model is able to deal with changes in the geometry, 
whereas POD is only able to deal with the changes covered by the snapshots.
This effect is shown in Section~\ref{sec:APOD}, more precisely, in Fig.~\ref{fig:l2_abs}. 
The POD model provides the displacement in every node of the cube, 
thus most quantities of interest which depend on the displacement are easily computed. 
In contrast, the HTA is only suitable to approximate a single quantity of interest. 
Consequently, another HTA has to be established, if another quantity of interest has to be investigated.
The advantage of this highly specialized approximation is the 
low cost of computation later on.
Construction times,
the time needed to evaluate $100$ random parameters as well as 
the convergence behaviour for the estimation of the mean and variance
are analyzed and compared in the following (see Tab.~\ref{tab:construction_and_run_times}). 

\begin{table}
	\centering
	\begin{tabular}{r|c|c}
		& constr. time [\SI{}{\s}] & $1000$ eval.  [\SI{}{\s}]\\
		\hline
     FEM  & -			&	32722.4	     \\	
     POD  & 130.89			&	19869.7	     \\	% Construction time = 32722.4(1000 FEMs)/1000 * 4(precalculations)	
    APOD  & 130.89			&	22320.8	     \\	% Construction time same as POD
     HTA of FEM  & 41497.5		&	0.245135	     \\	
     HTA of POD  & 25260.6		&	0.241672	     \\	
    HTA of APOD  & 28343.8		&	0.238908	     \\	

	\end{tabular}
	\caption{An overview of the computation times for construction and evaluation of (A)POD and HTA. With $8^3$ finite elements.}
	\label{tab:construction_and_run_times}
\end{table}

Naturally, the construction time of the unreduced FEM model is zero. 
The construction time of APOD and POD is relatively low, 
since only a handful of FE simulations are needed to generate the snapshots. 
For the construction of the HTA, around $10^3$ evaluations have to be carried out, 
therefore the construction time is very high in comparison, 
even if POD is used as basis. 
% It should be remembered, that the use of POD as underlying mapping results in a loss of accuracy as
% pointed out in Section~\ref{sec:approx_quality_hta}.
Note, that the times given in Tab.~\ref{tab:construction_and_run_times} are subject to fluctation, 
since the number of entries needed to construct a HTA depends on the random submatrices in 
the computation of the generalized cross approximation.
Further, it has to be taken into account that for theevaluation of a POD or APOD model an iterative solver with an adaptive number of iteration steps is used.
Yet, the orders of magnitude are consistent for all performed experiments.

\begin{figure}
	\centering
% 	\setlength\figureheight{0.35\textwidth*5/7}
% 	\setlength\figurewidth{0.35\textwidth}
% % 	\InputIfFileExists{Auswertungen/1000randoms/time_for_evaluations.tikz}{}{\textbf{!! Missing graphics !!}}
% 	\InputIfFileExists{time_for_evaluations_2500.tikz}{}{\textbf{!! Missing graphics !!}}
% 	\includegraphics[width=7cm]{pictures/hta/dspl_minmax_el_8.png}
     \includegraphics[width=7cm]{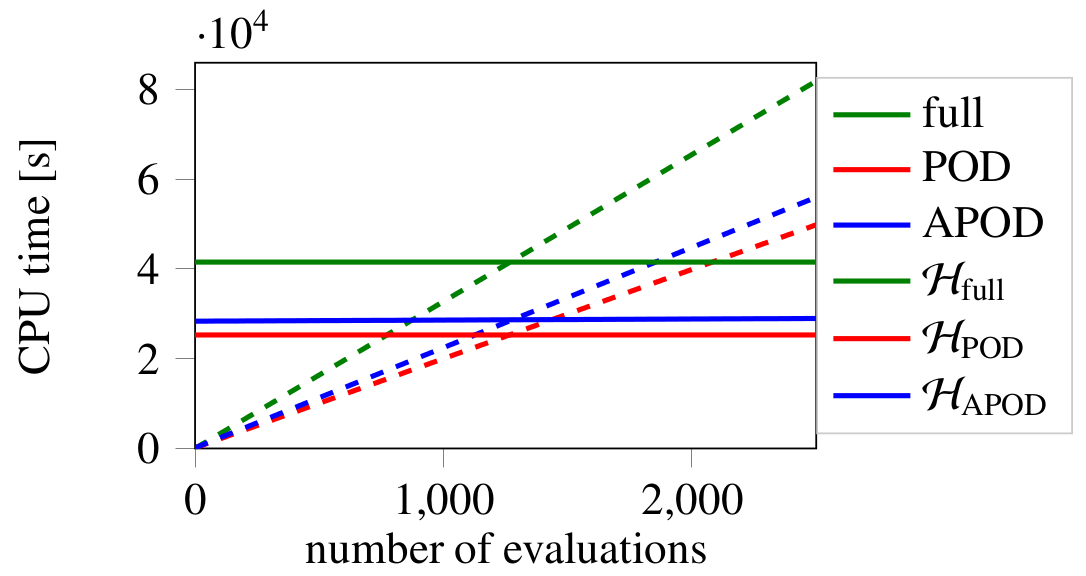}
	\caption{ Comparison of the computational time over the number of evaluations  }
	\label{fig:break_even_points}
\end{figure}

The high construction times of HTA  pay off, as it is seen for the cost of $100$ evaluations.
In Fig.~\ref{fig:break_even_points} the total CPU time is plotted over the number of evaluations.
For instance for 1000 evaluations, the use of HTA does not really pay off, 
because the construction time for establishing the $\mathcal{H}-$Tensor is very high in comparison to the CPU time needed to comput the results of 500 evaluations.
Already for 2500 evaluations, the picture has changed to the opposite.
Please note that a number of 2500 evaluations is easily reached and in the context of uncertainty quantification a very small number of evaluations.
% Keeping the differences in accuracy in mind, 
% one is able to choose the right model for the application at hand, 
% depending on the number of evaluations needed. 

\begin{figure}
	\centering
% 	\setlength\figureheight{0.35\textwidth*5/7}
% 	\setlength\figurewidth{0.35\textwidth}
% 	\InputIfFileExists{Means_over_samples.tikz}{}{\textbf{!! Missing graphics !!}}
% 	\InputIfFileExists{Auswertungen/1000randoms/Means_over_samples_mod2.tikz}{}{\textbf{!! Missing graphics !!}}
% 	\includegraphics[width=7cm]{pictures/hta/dspl_minmax_el_8.png}
     \includegraphics[width=7cm]{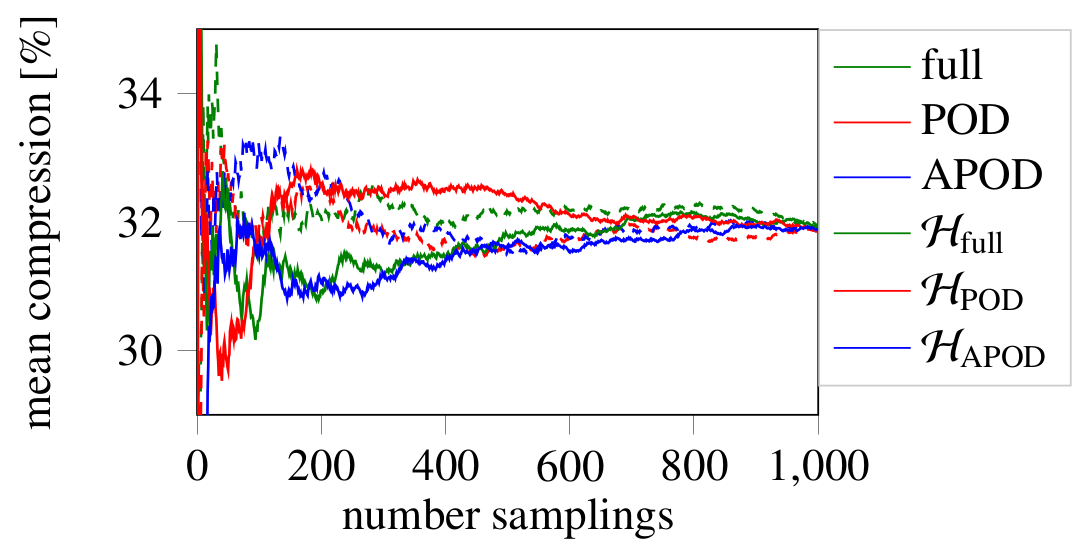}
	\caption{Convergence to the mean for all models up to $1000$ samplings.}
	\label{fig:means_over_N_rnd_samples}
\end{figure}

\begin{figure}
	\centering
% 	\setlength\figureheight{0.35\textwidth*5/7}
% 	\setlength\figurewidth{0.35\textwidth}
% 	\InputIfFileExists{variance_over_samples.tikz}{}{\textbf{!! Missing graphics !!}}
% 	\InputIfFileExists{Auswertungen/1000randoms/variance_over_samples_mod2.tikz}{}{\textbf{!! Missing graphics !!}}
% 	\includegraphics[width=7cm]{pictures/hta/dspl_minmax_el_8.png}
     \includegraphics[width=7cm]{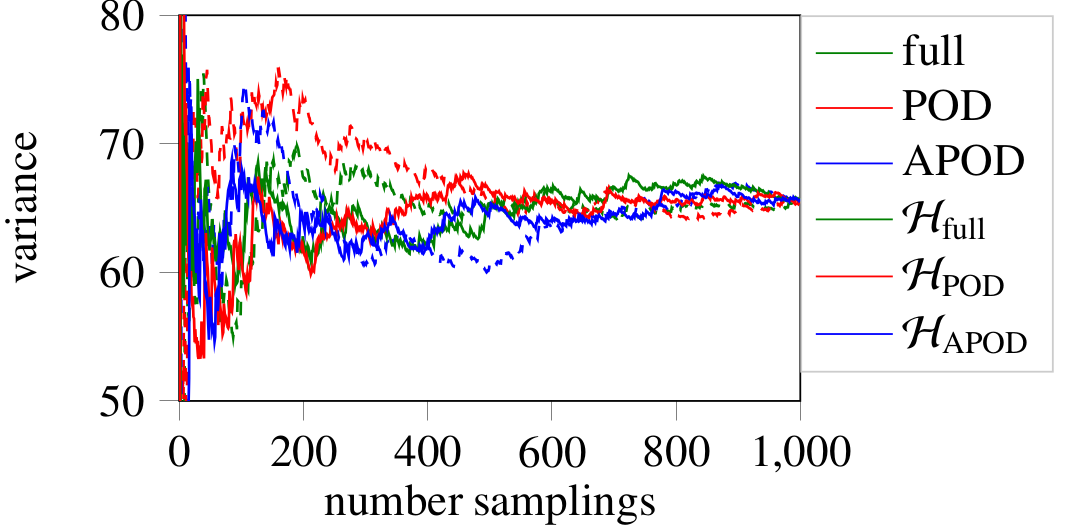}
	\caption{Convergence to the variance for all models up to $1000$ samplings.}
	\label{fig:variance_over_N_rnd_samples}
\end{figure}

The mean and the variance of the maximal displacement are shown in Fig.~\ref{fig:means_over_N_rnd_samples} and 
Fig.~\ref{fig:variance_over_N_rnd_samples}, respectively
For both, the mean and the variance, it can be observed that all models converge rapidly to a similar value.
% In Tab.~\ref{tab:phi_all} it is found that the values comply with the relative mean errors found in Tab.~\ref{tab:errors_all}.
The differences of the mean values in Tab.~\ref{tab:phi_all} are in the same magnitude as the errors found in Tab.~\ref{tab:errors_all}.

% This work is closed with a brief conclusion and an outlook in the following section.

\begin{table}
 \begin{tabular}{c|c|c|c|c}
  &       & std       &     & \\
  & mean  & deviation & max & min \\
  & [\%]  & [\%] & [\%] & [\%] \\
  \hline
$\varphi_{\mathrm{full}}$ &   31.936 & 8.1107 & 50.6728 & 20.6467 \\
$\varphi_{\mathrm{POD}}$ & 31.8818 & 8.0875 & 50.0193 & 20.6248 \\
$\varphi_{\mathrm{APOD}}$ & 31.8853 & 8.0926 & 50.1861 & 20.6199 \\
$\varphi_{\mathcal{H}_{\mathrm{full}}}$ & 31.936 & 8.1107 & 50.6727 & 20.6467 \\
$\varphi_{\mathcal{H}_{\mathrm{POD}}}$ & 31.8828 & 8.088 & 50.0455 & 20.6186 \\
$\varphi_{\mathcal{H}_{\mathrm{APOD}}}$ & 31.886 & 8.0967 & 50.1872 & 20.6196
\end{tabular} 
 \caption{QoI for different methods from $1000$ random parameters on the tensorized grid.}
 \label{tab:phi_all}
\end{table}

\section{Conclusion and Outlook}
\label{sec:outlook}
In this work three levels of modeling were presented. 
The original FE model which is here considered as full model was first reduced by 
(adaptive) proper orthogonal decomposition and then by low rank approximation in
hierarchical tucker format. 
It should be emphasized that the procedure is applicable to arbitrarily geometrically and materially non-linear problems. 
The present paper is, however, restricted to Neo Hookean behaviour, 
where uncertainties in geometric and material parameters were investigated.

First of all, both model reduction methods were assessed with respect to their approximation quality.
It was discovered that the accuracy of the A(POD) has to be improved by additional snapshots.
An efficient procedure to find these additional snapshots by means of a novel combination of a greedy algorithm with HTA was developed. 
The accuracy obtained from the automatically chosen snapshots was compared with the result of the manually chosen snapshots. 
Obviously, the method is promising, in particular for complex examples where the manual choice of the snapshots is very difficult.

Beside the approximation quality, the construction and evaluation times of (A)POD and HTA were compared.
The results indicate that using a low rank approximation of a (A)POD model is a viable combination of methods.
Since the error $\epsilon_{\mathrm{(A)POD},\mathcal{H}_{\mathrm{(A)POD}}}$ is several magnitudes smaller 
than the error $\epsilon_{\mathrm{full},\mathrm{(A)POD}}$, 
the approximation $\varphi_{\mathcal{H}_{\mathrm{(A)POD}}}$ 
becomes a useful surrogate for $\varphi_{\mathrm{(A)POD}}$ which in turn is a surrogate for the full model.
The relation between construction times of $\mathcal{H}_{\mathrm{(A)POD}}$ and $\mathcal{H}_{\mathrm{full}}$ is proportional to 
the relation of evaluation times of an (A)POD and a full simulation. 
Whereas the evaluation time of a low rank tensor is negligible in relation to one evaluation of the (A)POD model.
For large scale problems both the relation of construction and evaluation times will improve significantly in favour of this strategy.

In the development of this work, many ideas for future work emerged. 
These ideas were excluded here, in order to focus on the relation between POD and HTA.
With regard to modern methods in uncertainty quantification, like multi-level and quasi Monte Carlo methods, 
stochastic collocation and stochastic Galerkin methods,
one idea would be to incorporate our reduced models.
In terms of polymorphic uncertainty quantification, 
an attempt to incorporate fuzzy-stochastic numbers or other non-probabilistic methods could 
improve the representation of epistemic uncertainties.





% \bibliography{sample}
% \bibliographystyle{elsarticle-num-names}
\section{References}

\bibliographystyle{elsarticle-harv}
\bibliography{bib}

\end{document}